\documentstyle[amssymb]{article}
\begin{document}
%
\def\C{{\Bbb C}}
\def\F{{\Bbb F}}
\def\K{{\Bbb K}}
\def\N{{\Bbb N}}
\def\O{{\Bbb O}}
\def\P{{\Bbb P}}
\def\R{{\Bbb R}}
\def\s{{\Bbb S}}
\def\Q{{\Bbb Q}}
\def\Z{{\Bbb Z}}
\def\J{{\cal J}}
\def\B{{\frak B}}
\def\Im{{\rm Im}}
\def\Isom{{\rm Isom}}
\def\Exp{{\rm Exp}}
\def\D{{\rm D}}
\def\PSL{{\rm PSL}}
\def\PGL{{\rm PGL}}
\def\SL{{\rm SL}}
\def\Abs{{\rm Abs}}
\def\Aut{{\rm Aut}}
\def\Edge{{\rm Edge}}
\def\Vert{{\rm Vert}}
\def\SO{{\rm SO}}
\def\Sp{{\rm Sp}}
\def\SYMP{{\rm SYMP}}
\def\tr{{\rm tr}}
\def\SU{{\rm SU}}
\def\arctanh{{\rm arctanh}}
\def\GL{{\rm GL}}
\def\im{{\rm im}}
\def\ker{{\rm ker}}
\def\M{{{\cal M}}}
\def\H{{\cal H}}
{\bf \LARGE Notes on affine isometric actions of discrete
groups}

\vspace{22pt}

{\Large Yurii A. Neretin}

\vspace{22pt}

Moscow State Institute of Electronics and Mathematics, Bolshoi Triohsvyatitelskii, per. 3/12, Moscow-109028, Russia

Max-Planck-Institut f\"ur Mathematik, Gottfried-Claren-Str. 26, 53225 Bonn, Germany

\vspace{22pt}
neretin@mpim-bonn.mpg.de,

 neretin@main.mccme.rssi.ru
\vspace{22pt}
%
%
%
%
%
%
%
%


{\bf Contents}

{\it 1. Affine isometric actions of semisimple groups of rank 1.

2. General remarks on affine isometric actions.

3. Rigidity for affine isometric actions.

4. $\R$-trees and related affine isometric actions

5. General remarks (continuation).

6. Fock representation of semigroup of probabilistic measures on group

7. Group of automorphisms of bundle.

8. On representations of group of diffeomorphisms of Riemann surface

9. On central extensions of groups of automorphisms of the bundles.

10. Immobile sets and immobile functions.}

\newpage
\newcounter{yyy}
\newcounter{xxx}

\newcommand{\no}{\eqno{(\arabic{yyy}.\arabic{xxx})}\refstepcounter{xxx}}

\setcounter{xxx}{1}

These notes contain several more or less simple facts about affine isometric actions (see below subsection 0.1) of discrete groups. Mainly we are interested by lattices in rank 1
semisimple groups over $\R$ and over local fields. It is well-known, that
affine isometric action is a argument of several interesting functors.
Using affine isometric action it is possible to construct representations of
current groups (see [Ara], [Gui1], [VGG1-3]) of semigroups of probabilistic measures on a group (see[Ner2]), of groups
of diffeomorphisms(see [Ism4]), of categories of bistochastic measures
(see [Ner1-2]), and also central extensions of groups of some groups of diffeomorphisms and some other groups.

Unfortunately the majority of interesting groups have no affine actions at all.
The purpose of these notes is to show that discrete groups acting
on hyperbolic spaces have many natural affine isometric actions and to obtain
some corollaries from these constructions.

It is well-known that unitary representations of discrete groups is a dangerous subject and it contains many pathological phenomena. Neverless this topic contains some optimistic results, for instance unitary rigidity
theorems of Figa-Talamanca--Picardello type (see [FTP] and
  [BS], [CS]) or existence of holomorphic wavelets
(see below ss.3.1)\footnote{Of course there are some other
interesting phenomena, see [VDN]}. We try to continue this 'nonpathological line". We also show that in some cases "affinization"
of unitary representation allows to avoid some unitary upleasent phenomena. For instance it is possible to make
irreducible affine action
from factor-representation of discrete group.

{\bf 0.1. Affine isometric actions.}
Let $G$ be a topological group. Let $H$ be a (complex or real) Hilbert space. {\it A affine isometric action} (we will also use a term {\it AI-action})
$\Pi = (\pi,\gamma)$
of the group $G$ on $H$ is an action of $G$ by affine isometric transformations
of the space $H$, i.e. transformations having the form

$$\Pi (g) v= \pi (g) v + \gamma (g) \no\label{0.0}$$
where $v\in H$, operators
$\pi(g)$
are unitary operators and $\gamma(\cdot)$ is a function
$G \rightarrow H$. 
The equality

$$\Pi (g_1g_2) v = \Pi (g_1) \Pi (g_2) v$$
implies 

$$ \pi(g_1g_2) = \pi(g_1) \pi(g_2) \no\label{0.1} $$
$$\gamma(g_1 g_2) = \pi(g_1) \gamma (g_2) + \gamma(g_1) \no\label {0.3} $$
The condition (0.2) means that
$\pi(\cdot )$
is a unitary linear representation of $G$. Evidently
$\gamma (e) =0$,
hence (0.3) implies

$$\gamma (g^{-1}) =-\pi(g)^{-1} \gamma(g) \no\label {0.4}$$

{\bf 0.2.}
Let
$\pi$
be a unitary representation of $G$ in a Hilbert space $H$. Let
$w \in H$. Let us fix a vector $w\in H$. 
Then the formula

$$\Pi(g) v= \pi(g) v+ (\pi(g) w -w) \no\label {0.5}$$
defines an isometric action of $G$. The action (0.5) is not interesting: it becomes linear after the shift of the origin to the point $w$.

Now two following natural questions arise.

{\bf Question 0.1.}
Let
$\pi$
be a unitary representation of $G$. Does there exist a function
$\gamma : G \rightarrow H$
satisfying the condition (0.3) which is not representable in the form

$$\gamma (g) = \pi(g) w -w ? \no\label{ 0.6}$$
{\bf Question 0.2.} Let a hilbert space $L$ and a unitary representation
 $\pi$ of the group $G$ are given.
Is it possible to describe cohomology group
$H^1(G,L):=Z^1(G,L)/K^1(G,L)$ where the cochain group
$Z^1(G,L)$ consist of functions $G\rightarrow L$ satisfying the condition (0.3)
and the coboundary group $K^1(L,G)$ consist of all function having the form
(0.6)

{\bf 0.3}
The answer on Question 0.1 for the semisimple groups rarely is affirmative.

{\bf Theorem 0.3}{\it
a) Let a connected real semisimple group $G$ has a nontrivial affine isometric action. Then
$G = SO(1,n)$
or
$SU(1,n)$. \newline

b) Let a semisimple group $G$ over nonarchimedian local field has a nontrivial affine isometric action. Then rank of $G$ equals 1.}

{\bf Theorem 0.4.}{\it
a) Let $n>2$. Then the group
$SO(1,n)$
has only one nontrivial isometric affine action (in particular it is valid in the case
$\PSL_2(\C) = SO_0(1,3))$.\newline

b)The group
$SU(1,n)$
has only two nontrivial affine isometric actions (in particular it is valid for
$\PSL_2(\R) = SO_0(1,2))$.}

These theorems are well-known (and I don't know 
their history, see partial result in [Kaz],
 the proof of the first theorem see [Mar], III. 5, see also [Kaz];
 for real case see [Gui2];
for second theorem see [Gui2]. The construction of affine actions for Lie groups see below $\S$ \, 1.

%

{\bf 0.4.}
By this reason it seems that affine isometric actions is an exotic phenomena which is not very interesting. On other hand it turn out
that affine isometric actions are related to many algebraic, analytic and probabilistic phenomena and constructions: representations of current groups (Araki multiplicative integral, see [Ara], [VGG1] - [VGG3], [Ber], [Gui1], [PS], [Del1] - [Del2], [Ner2]), conditionally positive defined functions (see [Gui1], [PS], [VK], [Mar]), categories of $G$-stochastic measures ([Ner 1-2], representations of semigroup of probabilistic measures on groups (see [Ner2]), representations of groups in Pontryagin spaces, representations of various
 infinite-dimensional groups (see [Ism 1-2], [Ols 2-3], [Ner 1-2], [ATH]).

{\bf 0.5.}
The purpose of these notes is to show that discrete groups which act on hyperbolic spaces have many irreducible affine isometric representations. \newline

In $\S$ \, 1 we describe affine isometric representations of groups
$SO_0(1,n), SU(1,n)$, 
$\PSL_2(\Q_p)$ and pf the group of automorphisms
of Bruhat-Tits tree.
All these constructions are known. \newline

In $\S$ \, 2 we discuss general formal definitions and properties of affine isometric actions. \newline
In $\S$ \, 3 we prove the following rigidity theorem. Let $G$ be a semisimple group and
$\Gamma$ be a cocompact lattice. Then the restriction of each irreducible affine isometric action 
$\Pi$
of $G$ to
$\Gamma$
is irreducible. We also show that for
$G = \PSL_2(\R)$
the restriction
$\Pi|_{\Gamma}$
``remember'' embedding
$\Gamma \rightarrow \PSL_2(\R)$.
This theorem is a variation of Figa-Talamanca--Picardello--Bishop--Cowling--Steger theorems
(see [BS], [CS]) on unitary rigidity.

By [CM], [Mor] the space of homomorphisms of discrete group
$\Gamma$
to
$SO(1,n)$ (defined up to conjugation)
has natural compactification by the actions of
$\Gamma$
on
$\R$-trees (for $n=2$ and for the space 
of cocompact lattices
this compactification coincide with the Thurston compactification
 of Teichm\"uller space). \newline
In $\S$ \, 4 we show that for each action of
$\Gamma$
on
$\R$
-tree there exist a canonical irreducible isometric affine action of
$\Gamma$
in some Hilbert space associated to the tree.
Nonisomorphic actions of
$\Gamma$
on
$\R$
-trees generates nonisomorphic affine isometric actions. Notice that associated linear unitary representations of
$\Gamma$
are factor-representations and they ``forget'' action of
$\Gamma$
on tree. \newline

Our $\S$ \, 5 contains continuation of preliminaries.
We discuss relation of affine isometric actions with
infinitely divisible representations.

of \newline
In $\S$ \, 6 for each irreducible affine isometric action of lattice
$\Gamma$ we construct an irreducible representation of semigroup of probabilistic measures on
$\Gamma$
in the boson Fock space. \newline

In $\S$ \, 7  $\S$ 9 we discuss the group
$\B(G)$
of (measurable) automorphisms of principal $G$-bundle, we construct its representations (it is well-known) and central extensions (by strange way such extensions never wasn't discuss). \newline

In $\S$ \, 8 we discuss embeddings of groups of volume preserving diffeomorphisms to groups
$B(\Gamma)$
and associated representations. 

It is clear that it is impossible to describe all affine isometric actions of a given discrete groups
(of course if the set of these actions is not empty).
Neverless it seems that Question 0.2 can be nonpathological.
The discussion of this question is the subject of
$\S$ 9.

I am grateful to Gr. A. Margulis, R. Howe, V.V.Fock for 
discussions of the subject and  references.

\section{ Affine sometric actions of simple Lie groups of rank 1} 

\stepcounter{yyy}
\setcounter{xxx}{1}

Here we discuss constructions of affine isometric actions
of simple groups o rank 1.
  The most interesting cases for us
 are affine isometric actions of $\PSL_2(\R)$ (ss.1.1) and of  group
of automorphisms of Bruhat-Tits tree (ss 1.4).

{\bf 1.0. Preliminary remarks.}
Consider an affine isometric action
$\Pi = (\pi,\gamma)$

$$\Pi(g) v = \pi(g) v+ \gamma (g)$$
of a group $G$ in a Hilbert space $H$. Consider an operator

$$A(g) = \left(\begin{array}{cc} 1 & \gamma (g) \\ 0 & \pi (g) \end{array}\right)$$
in
$\C \oplus H$.
Condition
$\Pi(g_1 g_2) = \Pi(g_1) \pi(g_2)$
is equivalent to the condition

$$A(g_1g_2) = A(g_1) A(g_2)$$
Constructions of subsections 1.1 - 1.4, 1.6 below  have the following form. There exist some space
$\C \oplus H$
and linear nonunitary representation
$A(g)$
of the group $G$ in
$\C \oplus H$
such that subspace $H$ is
$A(g)$-invariant, operators
$A(g)$
are unitary in the subspace $H$ and the action of $G$ 
in the quotient space 
$(\C \oplus H)/H$
is trivial.
Consider a nonzero vector
$w \in \C \oplus 0$.
Then

$$A(g)w -w \in H$$
and hence the formula

$$\Pi_s(g) v = A(g) v + s(A(g)w-w)$$
defines an affine isometric action $\Pi_s$ of 
the group $G$ in the Hilbert space $H$.

In subsections 1.2, 1.4, 1.6 this space
$\C \oplus H$
is quite visible (for construction of subsections 1.1, 1.3 this space is easily visualized for dual representations).

{\bf 1.1. Isometric actions of
$\PSL_2(\R)$} (see [VGG1]) We realize the group
$\PSL_2(\R)$
as the group
$\PSL_2(\R) = \SU(1,1)$
of matrices

$$g = \left(\begin{array}{cc} a & b \\ \bar b & \bar a 
\end{array}\right) ;\qquad |a|^2 - |b|^2 =1$$
We will consider matrix elements

$$a =a (g) \qquad b =b(g)$$
as functions on the group
$\PSL_2(\R))$. We will wright
$$
g=\left(\begin{array}{cc}
a(g)&b(g)\\ {\overline a(g)} &{\overline b(g)}
\end{array}\right)
$$

The group
$\PSL_2(\R)$
acts in the disc
$D:|z| < 1$
by the Moebius transformations

$$z \longmapsto z^{[g]} := \frac {az+b} {\bar bz + \bar a}$$
Denote by $H$ the Bergman space in the disc
$|z| < 1$,
i.e. the space of holomorphic functions in disc
$|z| < 1$
equipped by
$L^2$
-scalar product

$$\langle f_1, f_2 \rangle = \frac 1 \pi \int_{|z|<1} f_1(z) \overline{f_2(z)} dz d\bar z$$
It is easy to see that
$$
\langle z^n,z^m\rangle=\frac{1}{n+1}\delta_{n,m}
$$
The group
$\PSL_2(\R)$
acts in the Hilbert space $H$ by affine isometric transformations

$$\Pi_s(g) f(z) = f(z^{[g]})(\bar b z+\bar a)^{-2} + s\cdot \gamma (g)$$
where
$s > 0$
and

$$\gamma (g) = \frac{\bar b} {\bar b z + \bar a} \no\label {1.1}$$
It is easy to see that

$$\langle \gamma \left(\begin{array}{cc} a & b \\ \bar b & \bar a \end{array}\right);
 \gamma \left(\begin{array}{cc} \alpha &\beta \\ \bar \beta & \bar \alpha \end{array}\right) \rangle
 = -\ln \left( 1 - \frac {\bar b \beta} {\bar a \alpha} 
\right)$$
It is useful to rewrite this formula in the form

$$\langle \gamma (g_1), \gamma (g_2) \rangle = -\ln \frac{\overline{a(g_1 g_2^{-1})}}
{\overline{a(g_1)} \overline { a(g_2^{-1})}} 
\no\label{ 1.2}$$

{\bf Remark 1.1.} 
Denote by  $u$ the expression under the logarithm. Then
$|u-1| < 1$
and hence this logarithm is a well-defined single-valued function. Nevertheless
$\ln(a(\cdot ))$
is not well-defined and we cannot express the formula (1.2) in the form

$$-\left[\ln \, \overline{a(g_1g^{-1}_2)} -\ln \, \overline{a(g_1)} - \ln \, \overline{a(g_2^{-1})}\right]$$

Supposing 
$g_1 = g_2$
we get

$$\left\| \gamma \left(\begin{array}{cc} a & b \\ \bar b & \bar a \end{array}\right) \right\|^2
 = -\ln(1-\left| \frac b a\right|^2)
 =2 \, \ln \, |a|$$
Denote by
$d(\cdot, \cdot)$
the standard $\PSL_2(\R)$-invariant Poincare metric on the disc
$|z| < 1$:

$$d(z_1, z_2) = {1\over 2} \ln
\frac{|1-{\overline z_1}z_2|+|z_2-z_1|}
{|1-{\overline z_1}z_2|-|z_2-z_1|}\no\label{distan} 
$$

Denote by $\delta(g)$ the distance between points $0$ and
$0^{[g]}$:

$$\delta (g) =d (0,0^{[g]}) = 
\arctanh\left( \left | \frac b a \right | \right)=
\ln(|a|+|b|) $$

Simple calculation shows

$$\| \gamma (g) \|^2 = 2\delta (g) -2 \ln 2 + o(1),\qquad
|a|\rightarrow\infty \no$$

{\bf 1.2. Isometric actions of the groups $\SO_0(1,n)$.} (see [VGG2])
The group
$\SO_0(1,n)$
is the group of real matrices
$\left(\begin{array}{cc} a & b \\ c & d \end{array}\right)$
of dimension
$(1 + n) \times (1+n)$
such that

$$ \left(\begin{array}{cc}
 a & b \\ c & d \end{array}\right) \left(\begin{array}{cc} -1 & 0 \\ 0 & 1 \end{array}\right) \left(\begin{array}{cc} a & b \\ c & d \end{array}\right)^t = \left(\begin{array}{cc} -1 & 0 \\ 0 & 1 \end{array}\right) $$
$$
a > 0;\qquad \det \left(\begin{array}{cc} a & b \\ c & d \end{array}\right) = +1 $$
Let
$S^{n-1}$
be the sphere
$x^2_1 + \ldots + x^2_n =1$
in
$\R^n$.
The group
$\SO_0(1,n)$
acts on
$S^{n-1}$
by the conformal transformations

$$x \longmapsto (a + xc)^{-1} (b + xd)$$
where
$x = (x_1 \ldots x_n) \in S^{n-1}$
is a row-matrix.

Denote by $Z$ the space of
$C^{\infty}$
-functions $f$ on
$S^{n-1}$
satisfying the condition

$$\int_{S^{n-1}} f(x) dx =0$$
Consider the positive defined scalar product

$$\langle f,g \rangle = -\int_{S^{n-1}} \int_{S^{n-1}} \ln \mid x-y \mid f(x) \overline{g(y)} dxdy$$
in the space $Z$. Denote by $H$ the completion of $Z$ with respect to this scalar product. The group
$\SO_0(1,n)$
acts in $H$ by linear unitary transformations

$$\pi \left(\begin{array}{cc}
 a & b \\ c & d \end{array}\right) f(x) = (a+xc)^{-(n-1)} f((a+xc)^{-1}(b+xd))$$.

{\bf Remark 1.2}
The representation
$\pi$
is the so-called limit of complementary series.

 The expression

$$\Pi_s\left(\begin{array}{cc}
 a & b \\ c & d \end{array}\right) f(x) = \pi \left(\begin{array}{cc}
a & b \\ c & d \end{array}\right) f(x) + s \left( (a+xc)^{-(n-1)} -1\right)$$
defines a affine isometric action of the group
$\SO_0(1,n)$
in $H$.

{\bf 1.3. Isometric actions of the groups $\SU(1,n)$} (see [Ber2], [VGG2])
The group
$\SU(1,n)$
is the group of complex 
$(1 +n)\times (1+n)$
-matrices
$\left(\begin{array}{cc}
a & b \\ c & d \end{array}\right)$
satisfying the conditions

$$\left(\begin{array}{cc}
 a & b \\ c & d \end{array}\right) \left(\begin{array}{cc}
 -1 & 0 \\ 0 & 1 \end{array}\right) \left(\begin{array}{cc}
a & b \\ c & d \end{array}\right)^* = \left(\begin{array}{cc}
 -1 & 0 \\ 0 & 1 \end{array}\right) $$
$$
\det \left(\begin{array}{cc}
 a & b \\ c & d \end{array}\right)
 =1 $$
Let
$B_n$
be the unit ball
$|z_1|^2 + \ldots + |z_n|^2 < 1$
in
$\C^n$.
The group
$\SU(1,n)$
acts on
$B_n$
by transformations

$$z \longmapsto z^{[g]} := (a + zc)^{-1} (b + zd)$$
where
$z = (z_1 \ldots z_n) \in B_n$
is a row-matrix.

Consider the matrix-valued
reproducing kernel (on reproducing kernels
see for instance [FK]).

$$K(\cdot , \cdot) : B_n \times B_n \rightarrow GL(n,\C)$$
defined by the formula

$$K(z,u) = \left( 1-\sum^n_{j=1} z_j \bar u_j \right)^{-1} \left(1- \left(\begin{array}{ccc} z_1  \bar u_1 &\ldots & z_n  \bar u_1 \\
\dots &\dots&\dots\\
z_1  \bar u_n & \ldots & z_n  \bar u_n \end{array}\right)
 \right)$$

$$
=\left( 1-\sum^n_{j=1} z_j \bar u_j \right)^{-1} 
\left( 1-\left(\begin{array}{c}(z_1\\ 
\vdots\\z_n\end{array}\right)
\left(\begin{array}{ccc}\bar u_1&\cdots&\bar 
u_2\end{array}\right)
\right)$$
This kernel defines some  Hilbert space $H$ of
$\C^n$-valued holomorphic functions on
$B_n $
(we consider elements of
$\C^n$
as vector-columns).

Consider the unitary representation
$\pi$
of the group
$\SU(1,n)$
in the space $H$ given by the formula

$$\pi \left(\begin{array}{cc} a & b \\ c & d \end{array}\right)
 f(z) = (a+zc)^{-1} (d-cz^{[g]})f(z^{[g]}).$$

In this formula $(a+zc)^{-1}$ is a scalar-valued function
and $(d-cz^{[g]})$ is a function with values in $n\times n$ matrices.

{\bf Remark 1.3}
The representation
$\pi$
is one of the highest weight representations of
$\SU(1,n)$.
This representation is degenerated, i.e. the space $H$ doesn't contain all polynomial functions
 $B_n\rightarrow \C^n$. The restriction of
$\pi$
to the group
$\SU(n)$
is the sum of the symmetric powers
$\bigoplus^{\infty}_{k =1} S^{k} R$
of standard representation $R$ of
$\SU(n)$
in
$\C^n$
(emphasize that
$k \ne 0)$.

Affine isometric action of
$\SU(1,n)$
in the hilbert space $H$ is given by the formula

$$\Pi_s \left(\begin{array}{cc} a & b \\ c & d \end{array}\right)
f(z) = \pi \left(\begin{array}{cc} a & b \\ c & d \end{array}\right)
 f(z) + s \cdot (a + zc)^{-1} \cdot c$$
where
$s > 0$. We remind that $c$ is a vector-column.

{\bf 1.4. Affine isometric actions of groups of automorphisms of Bruhat-Tits trees.}
Let $n$ be a positive integer. {\it Bruhat-Tits tree}
${\cal T}_n$
is the tree such that each vertex is incindent to
$(n+1)$
edges. Let
$s, s^{\prime}$
be vertices of
${\cal T}_n$.
The distance
$s,s^{\prime}$
is the length of (unique) chain joining $s$ and
$s^{\prime}$.

The {\it absolute} (the set of points on infinity)
$\Abs({\cal T}_n)$
 of the tree
${\cal T}_n$
is defined by an obvious way. Fix a vertex
$s^* \in {\cal T}_n$.
Then there is obvious one-to-one correspondence between points of
$\Abs({\cal T}_n)$
and infinite chains

$$s_0 =s^*, s_1, s_2, s_3, \ldots$$
(where
$s_j$
and
$s_{j+1}$
are incindent to the same edge and
$s_j \ne s_{j+2}$).
Define a metric
$\delta ( \cdot , \cdot )$
on
$\Abs({\cal T}_n)$.
Let
$x,y \in \Abs({\cal T}_n)$.
Let
$\ldots, t_{\alpha}, t_{\alpha + 1}, \ldots$
be a chain of vertices joining $x$ and $y$. Then

$$\delta (x,y) = n^{-\min(d(s^*,t_j))}$$
The space
$\Abs({\cal T}_n)$
with metric
$\delta$
is compact totally disconnected metric space.

 Define a {\it 
canonical measure}
$\mu$
on
$\Abs ({\cal T}_n)$
by the condition: measure of each ball with radius
$n^{-k}$
equals
$n^{-k}$.

{\bf Remark 1.4.} The metric $d(\cdot,\cdot)$ and the 
measure $\mu$ depends on the point
$s^*$,i.e. $d(\cdot,\cdot)=d_{s^*}(\cdot,\cdot))$ ,
$\mu=\mu_{s^*}$.

Denote by
$\Aut({\cal T}_n)$
the group of automorphisms of
${\cal T}_n$.
It is locally compact totally disconnected group. On representation theory of
$\Aut ({\cal T}_n)$
see [Ols 1], [FTN].

Let
$g \in \Aut({\cal T}_n), x \in \Abs ({\cal T}_n)$.
Define the {\it derivative}
$g^{\prime}(x)$.
Let
$s_0 =s^*, s_1, s_2, \ldots$
be the chain joining $s^*$ and $x$ and let
$t_0 =s^*, t_1, t_2,\dots$
be the chain joining $s^*$ and
$gx$.
Then there exists
$\alpha$
such that
$gs_j = t_{j+\alpha}$
for sufficiently large $j$. The derivative is defined
by the equality

$$g^{\prime}(x) = n^{\alpha}.$$

{\bf Remark 1.5.} Of course this is the Radon-Nikodim derivative of the measure
$\mu$ with respect to the transformation $g$.

Let
$H_0$
be the space of locally constant functions $f$ on
$\Abs({\cal T}_n)$
such that
$\int f d\mu =0$.
Consider the (positive defined) scalar product in
$H_0$
given by the formula

$$ \langle f_1,f_2 \rangle =-\int_{\Abs({\cal T}_n)} \int_{\Abs({\cal T}_n)} {\delta(x,y)}^{-1}f_1(x) 
{\overline{ f_2(y)}} d\mu (x) d\mu (y). $$
Let $H$ be the completion of
$H_0$
with respect to this scalar product. Isometric action of
$\Aut({\cal T}_n)$
in the Hilbert space $H$ is given by the formula

$$\Pi_s(g)f(x) = f(gx)g^{\prime}(x) + s(g^{\prime}(x)-1)$$

{\bf 1.5. On harmonic analysis on trees.}
Material of this section will be used only in subsection 4.8 and section 10.
Denote by $\Vert=\Vert({\cal T}_n)$ the set of vertices of ${\cal T}_n$,
and by $\Edge^*=\Edge^*({\cal T}_n)$ the set of oriented edges.
For each $\sigma\in\Edge^*$ denote by $\phi(\sigma)$
and $\psi(\sigma)$ the origin and the end of $\sigma$
respectively.

Consider spaces $l_2(\Vert)$ and $l_2^\pm(\Edge^*)$, elements of the latter space are $l_2$-functions on $\Edge^*$
satisfying the condition
$$f(-\sigma)=-f(\sigma)$$
where $(-\sigma)$ is the oriented edge inverse to $\sigma$.

Consider three $\Aut({\cal T}_n)$-intertwining operators.

1. {\it Laplace operator} $\Delta:l_2(\Vert)\rightarrow l_2(\Vert)$
is given by the formula
$$\Delta f(v):=\frac1p\sum_{d(v_j,v)=1} f(v_j)-f(v)$$

2. Operator $\bigtriangledown:l_2(\Vert)\rightarrow l_2^\pm(\Edge^*)$ is given by the formula
$$\bigtriangledown f(\sigma):=f(\psi(\sigma))-f(\phi(\sigma))$$

3. Operator $\blacktriangledown:l_2^\pm(\Edge^*)\rightarrow
l_2(\Vert)$ is given by the formula 
$$\blacktriangledown h(v)=
\sum_{\sigma_j:\phi(\sigma_j):=v} h(\sigma_j)$$

Obviously all operators $\Delta,\bigtriangledown ,\blacktriangledown$ are bounded and
$$\blacktriangledown\bigtriangledown =p\cdot\Delta ;\qquad 
\blacktriangledown^*=\bigtriangledown \no\label{inve}$$

Spectrum of Laplace operator in the space 
$l_2(\Vert)$ is known (see [Ols1],[FTN])
and 0 is not a point of spectrum, i.e $\Delta$ is
invertible. Hence (1.6) implies that image 
$\im(\bigtriangledown)$ of
$\bigtriangledown$ is a closed subspace and its
orthocomplement is $\ker(\blacktriangledown)$.

Consider the space $H$ defined in preceding subsection.
We will describe the canonical unitary operator
${{\cal P}}:H\rightarrow\ker(\blacktriangledown)$
 (Poisson transform) which intertwines (linear) actions
of $\Aut({\cal T}_n)$ in the spaces $H$ and $\ker(\blacktriangledown)$. Consider a edge $\sigma$.
Then the space $\Abs({\cal T}_n)$ is divided by natural way
on
two pieces $\Abs_\sigma^+$ and $\Abs_\sigma^-$, the points of
$\Abs_\sigma^+$ are more close to origin $\phi(\sigma)$ 
of $\sigma$ and
points of $\Abs_\sigma^-$ are more close to the end 
$\psi(\sigma)$.
Then
$${{\cal P}}f(\sigma):=C\cdot\left(\int_{\Abs_\sigma^+}f(x)d\mu_{\phi(\sigma)}-
\int_{\Abs_\sigma^-}f(x)d\mu_{\psi(\sigma)}\right)$$
where $C$ is some constant.

{\bf 1.6. Isometric actions of groups
$\SL_2$ over local fields}. Let
$\K$ be a locally compact nonarchimedian field. Let
$\cal O \subset \K$
be the ring of integers, let $I$ be the nonzero prime ideal in
$\cal O$.
Let
$\F = {\cal O}/I$
be the residue field, let
$q=p^{k}$
be the number of elements of $\F$.

{\bf Remark 1.6.} Let $\K=\Q_p$ be $p$-adic field. Then
 ${\cal O}={\Z}_p$ is the 
ring of $p$-adic integers and $I=p\cdot Z_p$.

A {\it lattice} in the linear space
$\K^2$
is
$\cal O$
-submodule generated by two linearly independent vectors. Two lattices
$\Lambda$
and 
$\Lambda^{\prime}$
are {\it equivalent} if there exists
$x \in \K$
such that
$\Lambda^{\prime} = x \cdot \Lambda$.
Now we will construct a graph
${\cal T}$.
The vertices of
${\cal T}$
are enumerated by equivalence classes of lattices. Vertices $s$ and
$s^{\prime}$
are adjacent if there exist representative lattices 
$\Lambda$
and
$\Lambda^{\prime}$
such that index of
$\Lambda^{\prime}$
in
$\Lambda$
is $q$. Then the graph
${\cal T}$
is Bruhat-Tits tree
${\cal T}_q$.
The space
$\Abs({\cal T}_q)$
is identified with projective line
$\P \K^1$
over
$\K$
(see [Ser2], [FTN]).

By construction the group
$\PSL_2(\K)$
acts on the graph
${\cal T} = {\cal T}_q$
and hence we obtained embedding

$$\PSL_2(\K) \rightarrow \Aut ({\cal T}_{p^k})$$
Thus we can restrict the affine isometric action of
$\Aut({\cal T}_{p^k})$
to
$\PSL_2(\K)$.

Now we will give another description of the same action. Let
$H_0$
be the space of  locally constant functions
$f : \K \rightarrow \C$
such that
$\int f=0$.
Consider the scalar product in
$H_0$
given by the formula

$$\langle f_1,f_2 \rangle = -\int_{\K} \int_{\K} \ln|x-y|f(x) \overline{f(y)} dxdy$$
Let $H$ be the completion of
$H_0$
with respect to this scalar product. Isometric action of
$\PSL_2(\K)$
in Hilbert space is given by the formula

$$\Pi_s\left(\begin{array}{cc} a & b \\ c & d \end{array}\right)
 f(x) = f \left( \frac {ax+b} {cx+d} \right) |cx+d|^{-2} + s \cdot \biggl( |cx+d|^{-2}-1 \biggr) .$$

\section{ General remarks on affine isometric actions}

\stepcounter{yyy}
\setcounter{xxx}{1}

{\bf 2.1. Equivalent isometric actions.}
Consider two affine isometric actions of a group $G$:

$$ \Pi_1(g) v = \pi_1(g) v + \gamma_1(g) $$
$$\Pi_2(g) w = \pi_2(g) w + \gamma_2(g) $$
of a group $G$ in Hilbert spaces
$H_1$
and
$H_2$.
The actions
$\Pi_1$
and
$\Pi_2$
are called {\it equivalent} if there exists a map

$$Rh = Th +q$$
where
$T : H_1 \rightarrow H_2$
is a unitary operator and
$q \in H_2$
such that

$$R \Pi_1(g) = \Pi_2(g) R$$
for each
$g \in G$.

{\bf Lemma 2.1.}{\it
Let the actions
$\Pi_1$
and
$\Pi_2$
be equivalent. Then

 a)  Unitary representations
$\pi_1$
and
$\pi_2$
of the group $G$ are equivalent. 

 b) The function

$$\delta (g) = \parallel \gamma_1(g) \parallel^2 - \parallel \gamma_2(g)\parallel^2$$
is bounded.}

{\bf Proof}
a) is obvious.

{\bf Proposition 2.2.}{\it
Let
$\Pi = (\pi,\gamma)$
be an affine isometric action of a group $G$. The following conditions are equivalent

 a) The action
$\Pi$
is equivalent to a linear action.

 b)  The function
$\parallel \gamma(g)\parallel^2$
is bounded on $G$.}

{\bf Proof.}
a) $\Rightarrow$ b) is a partial case of b) $\Rightarrow$ a)
in lemma 2.1. Consider the space
$\C \oplus H$
and linear operators in
$\C \oplus H$
given by the formula

$$A(g) = \left(\begin{array}{cc} 1 & \gamma (g) \\
0 & \pi (g) \end{array}\right)
$$
Consider the closed convex hull $K$ of vectors
$(1,\gamma (g)) \in \C \oplus H$.
The set $K$ is bounded and hence (see [Rud], 3.15) $K$ is weakly compact. The group of operators
$A(g)$
is equicontinuous. By Kakutani theorem (see [Rud], 5.11) there exist a fixed point
$(0,h) \in K$
for the group $G$ of transformations
$A(g)$.
Then
$h \in H$
is a fixed point for transformations
$\Pi (g)$.
This proves the proposition.$\blacksquare$

We say that an affine isometric action is {\it nontrivial}
if it is not equivalent to a linear action.

{\bf Corollary 2.3.}{\it
All isometric actions constructed in $\S$ \, 1 are not equivalent to linear representations.}

{\bf 2.2 Minimal subspaces.}
Let $H$ be a Hilbert space. By definition {\it affine subspaces} of $H$ are shifts of closed linear subspaces.

An affine isometric action
$\Pi$
of a group $G$ in $H$ is called {\it AI-reducible
(affine isometric reducible)} if there exist $G$-invariant affine subspace
$L \subset H$.
The restrictions of isometric actions to $G$-invariant affine subspaces are called {\it AI-subrepresentations} of
$\Pi$.

{\bf Remark 2.4.}
Each linear representation is AI-reducible. Indeed the origin is an invariant subspace.

{\bf Proposition 2.5.}{\it
Let a nontrivial isometric action
$\Pi$
has a irreducible AI-subrepresentation $L$. Then each AI-subrepresentation of
$\Pi$
contains $L$.}

{\bf Proof.}
Let
$M \ne L$
be an invariant affine subspace. Suppose
$M \cap L = \varnothing$.
Let
$l \in L$
be the nearest to $M$ point of $L$. Then $m$ is a $G$-stable point and hence the action 
$\Pi$
is equivalent to linear action.

Assume $M\cap L\ne\varnothing$. Then
$M \cap L$
be a $G$-invariant affine subspace of $L$. Bearing in mind that $L$ is irreducible  we obtain
$M \supset L$.$\blacksquare$

{\bf Example 2.6.}
Consider isometric transformation

$$Af(x) = e^{ix}f(x) +1$$
in
$L^2[-1,1]$.
Consider an affine subspace
$V_{\varepsilon} \subset L^2[-1,1]$
which consist of functions satisfying the condition

$$f(x) = \frac 1 {1-e^{ix}} \qquad \mbox{for} \qquad |x| > \varepsilon$$
Then
$V_{\varepsilon}$
is an AI-invariant affine subspace and
$\bigcap_{\varepsilon > 0} V_{\varepsilon} = \varnothing$
(the function 
$1/(1-e^{ix})$
is not an element of
$L^2$).
Hence the action $A$ has no minimal invariant subspace.

{\bf 2.3. Semigroup of probabilistic measures.}
Denote by
${\M}_0(G)$
the semigroup of probabilistic measures on $G$ with compact support.

Let
$\Pi = (\pi, \gamma)$
be an isometric action of $G$. We extend the action
$\Pi$
 to the action
$\tilde \Pi$
of
$\M_0(G)$
given by the formula

$$\tilde \Pi (\mu)v = \int_G \pi(g) v \, d\mu(g) + \int_G \gamma (g) d\mu (g)$$

Let
$\pi$
be a unitary representation of a group $G$ in a Hilbert space $H$. By definition a representation
$\pi$
{\it weakly contains} trivial representation if there exist vectors
$v_j \in H$
such that the sequence of functions

$$\alpha_j(g) = \langle \pi (g) v_j, v_j \rangle$$
converges uniformly on compacts to 1 (see [Dix], 18.1).

{\bf Proposition 2.7.}{\it
(see [Mar], III.1.3) The following conditions are equivalent

 i) The representation
$\pi$
doesn't contain  weakly trivial representation.

 ii) There exist
$\mu \in \M_0(G)$
such that the norm of the operator

$$\tilde \pi(\mu) = \int_G \pi(g) d\mu(g)$$
is less than 1.
}

{\bf Proposition 2.8.}{\it
(see [Dix], 18.3.6, [Mar], III.4) The following conditions are equivalent

 i) A group $G$ is amenable

ii)  Regular representation of $G$  contains weakly trivial representation.
}

{\bf Proposition 2.9.}{\it
Let
$\Pi = (\pi, \gamma)$
be an affine isometric action of a group $G$. Let the representation $\pi$
doesn't contain  weakly trivial representation. Then
$\Pi$
contains an
AI
-irreducible subrepresentation.}

{\bf Proof.}
Consider
$\mu \in \M_0(G)$
such that
$\parallel \tilde \pi(\mu)\parallel < 1$
(see proposition 2.7). Then
$\tilde \Pi(\mu)$
is contractive map. Let $w$ be the unique fixed point of
$\tilde \Pi (\mu)$.
Consider the minimal affine subspace $L$ containing all points
$\tilde \Pi(g) w$ (where $g \in G)$.
Then $L$ is AI-irreducible subrepresentation.$\blacksquare$

{\bf 2.4 Direct sums.}
Let
$\Pi_1 = (\pi_{1},\gamma_1)$
and
$\Pi_2 = (\pi_2,\gamma_2)$
be isometric actions of the group $G$ in the spaces 
$H_1$
and 
$H_2$.
The {\it direct sum}
$\Pi_1 \oplus \Pi_2$
is affine isometric action of the group $G$ on the space
$H_1 \oplus H_2$
given by the formula

$$g: (v,w) \mapsto (\pi_1(g)v + \gamma_1(g), \pi_2(g)v + \gamma_2(g))$$
where
$v \in H_1, w \in H_2$.

{\bf Remark 2.10.}
Usually subspaces
$H_1, H_2$
in
$H_1 \oplus H_2$
are not invariant.

{\bf Example 2.11.}
Let
$H_1 = H_2, \Pi_1 = \Pi_2$.
Then the diagonal
$\Delta \subset H \oplus H$
is $G$-invariant affine subspace (diagonal 
$\Delta$
consists of vectors
$(v,v) \in H \oplus H)$.

Let
$\Pi = (\pi, \gamma)$
be an affine isometric action of $G$ in a Hilbert space $H$. Let a linear representation 
$\pi$
be reducible, let
$H = H_1 \oplus H_2$
be  a decomposition into a sum of subspaces which are invariant with respect to operators
$\pi(g)$.
Let
$\pi_1, \pi_2$
be the restrictions of representation
$\pi$
on
$H_1, H_2$.
Let
$\gamma_1(g), \gamma_2(g)$
be the projections of
$\gamma$
to
$H_1, H_2$.
Then
$\Pi = (\pi, \gamma)$
is equivalent to the direct sum of actions
$\Pi_1 = (\pi_1,\gamma_1), \Pi_2 = (\pi_2,\gamma_2)$.

{\bf Example 2.12.}
Let
$G = \Z$.
Denote by
$\Pi_s$
the action of
$\Z$
on 
$\R$
given by the formula

$$\Pi_s(n)x = x+ s n$$
$(x \in \R, n \in \Z)$.
Denote by
$\Pi^{(k)}_{\lambda}$
action of
$\Z$
on
$\R^k$
given by the formula

$$\Pi^{(k)}_{\lambda}(n)(x_1,\ldots, x_k) = (x_1 + \lambda, \ldots, x_k)$$
$(\lambda \in \R, n \in \Z, (x_1, \ldots, x_k) \in \R^k)$.
Then

$$\Pi_{s_1} \oplus \ldots \oplus \Pi_{s_k} \simeq \Pi^{(k)}_{\sqrt{s^2_1 + \ldots + s^2_k}}$$
This example demonstrated that decomposition of
an affine  isometric action onto direct sum of irreducible actions is not unique.

{\bf Remark 2.13.}
Let
$\Pi$
be an affine isometric action of a group $G$ in Hilbert space $H$. Let $K$ be a $G$-invariant subspace. Let
$0 \in K$
(otherwise we will shift origin). Let $M$ be the orthogonal complement to $H$. Decompose our action to direct sums of actions
$(\Pi_1,\pi)$
and
$(\Pi_2,\pi_2)$
in $K$ and $M$ (relatively).
Then action
$(\Pi_2,\pi_2)$
of group $G$ is linear.

{\bf 2.5 Real and complex action.}
There exist some difference between affine isometric actions in real and complex Hilbert space (see $\S$ 6-9).

Let
$\Pi = (\pi,\gamma)$
be an affine isometric action of the group $G$ on a Hilbert space $H$.

We say that the {\it action is real} in two cases:

 a) If the Hilbert space $H$ is real.

 b) If the Hilbert space $H$ is complex and there exist
$\Pi (g)$
-invariant real subspace
$H_{\R} \subset H$
such that
$H = H_{\R} \oplus i \, H_{\R}$.

We say that {\it action is complex} if it is not real.

In $\S$ \, 1 we constructed one real action for the groups
$\SO(1,n)$, $\Aut({\cal T}_n)$, $\PGL(2,\K)$ (up to dilatations).
For the groups
$\SU(1,n)$
(in particular for
$\SU(1,1) \simeq \PSL_2(\R))$
we have two nondecomposable complex actions (in 
spaces of holomorphic and antiholomorphic functions). Realification of these two actions coincides and this is the unique real nondecomposable real action for
$\SU(1,n)$.
In the case
$n=1$ 
this action coincides with real actions
$\SO_0(1,2)$
mentioned several rows above.

 \section{ Rigidity for affine isometric actions} 

\stepcounter{yyy}
\setcounter{xxx}{1}

Let a Lie group $G$ and a lattice $\Gamma\subset G$ are given.
It is well-known that there exist many cases when
the group $G$ and the lattice $\Gamma$ have closely
related behavior, see for instance [Mar]. Theorems formulated below
(ss.3.1-3.2) are very simple representatives 
of the rigidity theorems of this type.

{\bf 3.1. Rigidity for unitary representations.}
We remind some known theorems. Let $\pi$ be an irreducible unitary representation of a Lie group $G$. Let
$\Gamma \subset G$
be a {\it lattice}. Recall that a {\it lattice}
$\Gamma \subset G$
is a discrete subgroup such that the volume of
$G/\Gamma$
is finite. A lattice is called {\it cocompact} if the space
$G/\Gamma$
is compact.

Let $\pi$ be a irreducible unitary representation of the group
G. Consider restriction $\pi{\big |}_\Gamma$ 
of  representation $\pi$
to the lattice $\Gamma$. Remind that representation
$\pi$ is named by {\it representation of discrete series}
if its matrix elements are in $L^2$.

{\bf Theorem 3.1.}
{\it  Let $\pi$ be not a representation of discrete series. Then the restriction
$\pi{\big |}_\Gamma$
of $\pi$ to
$\Gamma$
is irreducible. If
$\pi_1 \ne \pi_2$
are not elements of discrete series then
$\pi_1{\big |}_\Gamma$
and
$\pi_2{\big |}_\Gamma$
are not equivalent.}

See [BS], [SC]; first theorem of this type was obtained
in [FTP].

{\bf Theorem 3.2.}
([BS]). {\it Let
$G = \PSL_2(\R)$,
let $\pi_1$, $\pi_2$ be not  representation of discrete series. Let
$\Gamma$
be a discrete group and
$\psi_1 : \Gamma \rightarrow \PSL_2(\R), \psi_2 : \Gamma \rightarrow \PSL_2(\R)$
be embeddings such that
$\psi_1(\Gamma)$
and
$\psi_2(\Gamma)$
are lattices. Then representations
$\pi_1\circ \psi_1$
and
$\pi_2\circ \psi_2$
are equivalent if and only if $\pi_1\simeq \pi_2$ and embeddings
$\psi_1, \psi_2$
are conjugate.}

{\bf Theorem 3.3.}
([AS], [GHJ]){\it

 a) Let $\pi$ be an irreducible representation of discrete series. Then the restriction
$\pi{\big |}_\Gamma$
is a factor-representation of
$\Gamma$.

 b) Let $\pi$ be a representation of discrete series and 
$\mbox{\rm Dim}(\rho) $
be the formal dimension of $\pi$. Let all conjugacy classes of
$\Gamma$
be infinite. Then the restriction
$\pi{\big |}_\Gamma$
is a factor-representation which is quasi equivalent to the regular representation of
$\Gamma$
in
$l_2(\Gamma)$.
The dimension of
$\pi{\big |}_\Gamma$
over
$l_2(\Gamma)$
equals

$$\mbox{\rm Dim} (\pi) \cdot \mbox{\rm volume}(G/\Gamma) \no\label {3.1}$$}

{\bf Example 3.4}, see [GHJ].
Let
$G = \PSL_2(\R)$.
Consider the representations
$\pi_k$
of discrete series of the group
$\PSL_2(\R)$.
The space
$H_k$
of representation
$\pi_k$
is the space of holomorphic functions on the disk
$|z| < 1$ equipped
with scalar product

$$\langle f_1,f_2\rangle = \int \int f_1(z) \overline{f_2(z)} (1-|z|^2)^{2k-2} dzd\bar z$$
The operators
$\pi_k(g)$
are given by the formula

$$\pi_k(g)f(z) = f \left( \frac{az+b} {\bar b z+ \bar a} \right) (\bar b z+ \bar a)^{-2k}$$
Consider 
$\PSL_2(\R)$
-invariant measure on disk
$|z| < 1$
such that the square of a triangle is
$(\pi$
-(sum of angles)). Then dimension of
$\pi{\big |}_\Gamma$
over
$l_2(\Gamma)$
equals to

$$\lambda = \frac{2k -1} 4 \cdot \mbox{\rm square} \, (\PSL_2(\R)/\Gamma) \no\label{ 3.2}$$

{\bf Problem 3.5 (holomorphic wavelets).}
Let number
$\lambda$
be given by the formula (3.2) be integer. Then 
the representation
$\pi{\big |}_\Gamma$
is isomorphic to the sum of
$\lambda$
copies of
$l_2(\Gamma)$.
Hence there exist an orthonormal basis
$e_{\alpha}$
in the space $H_k$ of the unitary representation $\pi$ such that for each
$h_i \in \Gamma$,
$$\pi_k(h_i)e_{\alpha} = e_{\beta}$$
for some
$\beta$. Is it possible to construct such basis explicitly?

{\bf 3.2. Rigidity for affine isometric actions.}

{\bf Theorem 3.6.}
{\it Let $G$ be a locally compact topological group and the Haar measure on $G$ is both-side invariant. Let
$\Gamma$
be a cocompact lattice. Let
$\Pi = (\pi,\gamma)$
be an AI
-irreducible isometric action of $G$. Then the restriction of
$\Pi$
to
$\Gamma$
is AI-irreducible.}

{\bf Corollary 3.7.}{\it
Let
$\Pi = (\pi,\gamma)$
be an irreducible affine action of a discrete group
$\Gamma$.
Let
$\Gamma^{\prime} \subset \Gamma$
be a subgroup of finite index. Then the restriction 
$\Pi'{\big |}_
\Gamma$
is AI-irreducible action.}

{\bf Example 3.8.}
Let
$G = \PSL_2(\R), \PSL_2(\K)$
or
$\Aut({\cal T}_n)$ and $\Gamma\subset G$ be a cocompact lattice.
Consider the affine isometric action $\Pi_s = (\pi,s\cdot\gamma)$
of the group G described in
$\S$ \, 1. Then the restriction of the representation
$\pi$
to the lattice
$\Gamma$
is a factor-representation, nevertheless the restriction of
$\Pi_s$
to
$\Gamma$
is AI
-irreducible.

{\bf Proposition 3.9.}{\it
Let
$G = \PSL_2(\R)$
and let
$\Gamma = \Gamma_n$
be the fundamental group of the sphere with $n$ handles
$(n > 1)$.
Let
$\psi_1, \psi_2$
be  embeddings
$\Gamma_n \rightarrow \PSL_2(\R)$. Let $\psi_1(\Gamma),\psi_2(\Gamma)$ be lattices in $\PSL_2(\R)$
Let
$\Pi_s = (\pi,s\cdot\gamma)$
be the affine isometric action of 
$\PSL_2(\R)$
described in subsection 1.1($s>0$). Then

$$\Pi_s \circ \psi_1 \simeq \Pi_{s^{\prime}} \circ \psi_2$$
only in the case when

$s = s^{\prime}$ and embeddings $\psi_1, \psi_2$ are conjugate.}

{\bf 3.3. Induced isometric actions.}
Let $G$ be a topological group. Let
$\Gamma$
be a cocompact lattice. Let
$\Xi = (\xi,\lambda)$ be an affine isometric
 isometric action of
$\Gamma$
in a Hilbert space $V$.

Consider the space
$\cal L$
consisting of all measurable functions
$G \rightarrow V$
satisfying the condition

$$f(hg) =\Xi (h) f(g) \no\label{ 3.3}$$
for each
$h \in \Gamma, g \in G$.

Let
$\Omega \subset G$
be a compact fundamental domain for
$\Gamma$
(i.e.
$\Omega$
is compact and a natural map
$\Omega \rightarrow \Gamma \setminus G$
is bijective (up to the set of zero measure) measure preserving map). Consider the space $H$ consisting of all functions
$f \in \cal L$
satisfying the condition

$$\int_{\Omega} \parallel f(g) \parallel^2 dg < \infty. \no\label{ 3.4}$$
The space $H$ is a Hilbert space with respect to the scalar product

$$\langle f_1, f_2 \rangle = \int_{\Omega} f_1 (g) \overline{f_2(g)} dg$$
Consider the action of $G$ on
$\cal L$
given by the formula

$$p : f(g) \mapsto f(gp^{-1}) \no\label{ 3.5}$$
where
$p \in G$.

{\bf Lemma 3.10.}{\it 
The space $H$ is invariant with respect to the transformations (3.5).}

{\bf Proof.}
Consider all elements
$h_1, \ldots, h_{\alpha} \in \Gamma$
such that
$\Omega p^{-1} \cap h_j \Omega \ne \varnothing$.

$$ \int_{\Omega} \parallel f(gp^{-1}) \parallel^2 dg = 
\int_{\Omega p^{-1}} \parallel f(g) \parallel^2 dg = $$
$$= \sum_j \int_{\Omega p^{-1} \cap h_j \Omega} \parallel f(g) \parallel^2 dg = 
\sum_j \int_{\Omega \cap h^{-1}_j \Omega p^{-1}} \parallel f (h^{-1}_j g) \parallel^2 dg = $$
$$= \sum_j \int_{\Omega \cap h^{-1}_j \Omega p^{-1}} \parallel \xi (h^{-1}) f(g) +\lambda (h^{-1}_j) \parallel^2 dg =$$
$$= \sum^{\alpha}_{j=1} \int_{\Omega \cap h^{-1}_j \Omega p^{_1}} (\parallel f(g) \parallel^2 + 2Re \langle \xi (h) f(g), \lambda (h^{-1}_j)\rangle + \parallel \lambda (h^{-1}_j)\parallel^2)dg < \infty .\blacksquare$$

Denote by
$\Phi (p)$
the restriction of transformations (3.5) to $H$. The space $H$ is an affine subspace of 
$\cal L$.
Hence transformations
$\Phi (p)$
are affine transformations. A linear part of an affine transformation
$\Phi (p)$
is the representation of $G$ induced from representation
$\xi$
of 
$\Gamma$.
Hence
$\Phi (p)$
is an affine isometric action.

{\bf 3.4. Proof of Theorem 3.6.}
Consider an AI-irreducible affine isometric action
$\Pi = (\pi,\lambda)$
of the group $G$ in a Hilbert space $K$. Denote by
$\Xi$
the restriction of the action
$\Pi$
to the lattice
$\Gamma$.
Consider the isometric action
$\Phi$
of $G$ induced from the isometric action
$\Xi$.
Let $H$ be a space of the action
$\Phi$, i.e. the space constructed in ss.3.3..

We will construct a natural $G$-equivariant embedding
$\delta : K \rightarrow H$.
Let
$v \in K$.
Associate to
$v$
the function
$\delta_{v} : G \rightarrow K$
given by the formula

$$\delta_{v} (g) = \Pi (g) v$$
Evidently the function
$\delta_v$
satisfies the condition (3.3). Test the conditions (3.4):

$$ \int_{\Omega} \parallel \Pi (g) v \parallel^2dg = \int_{\Omega} \parallel \pi (g) v + \gamma (g) \parallel^2 dg = $$
$$= \int_{\Omega}(\parallel v \parallel^2 + 2 \, Re \, \langle \pi(g)v, \gamma (g) \rangle + \parallel \gamma (g) \parallel^2)dg $$
The integrand is continuous and the domain
$\Omega$
is compact. Hence the integral converges and so
$\delta_{\sigma} \in H$.

Thus we construct
AI-irreducible subspace
$\delta (K)$
in $H$.

Now suppose that the affine isometric action
$\Xi$
of the lattice
$\Gamma$
is reducible. Let
$L \subset K$
be a 
$\Gamma$
-invariant subspace
$(L \ne K)$.
Consider the isometric action
$\Psi$
of $G$ induced from the isometric action of $\Gamma$ in $L$. Evidently
$\Psi$
is a AI-subrepresentation in
$\Phi$.
Evidently
$\Psi$
doesn't contain the whole space
$\delta (K)$.
But
$\delta (K)$
is irreducible and this contradicts to Proposition 2.5.
$\blacksquare$

{\bf 3.5. Length function.}
Consider the disc
$D : |z| < 1$
and the distance
$d(\cdot , \cdot)$
given by (1.3) on $D$. Let
$g = \left(\begin{array}{cc} a & b \\ \bar b & \bar a 
\end{array}\right) \in SU (1,1) \simeq \PSL_2(\R)$.
Define the {\it length function}
$\ell(g)$
on
$\SU(1,1)$
by

$$\ell(g) =  {\mbox{\rm inf}}_{|z|<1}  d(z,z^{[g]})$$
Let
$\lambda, \lambda^{-1}$
be eigenvalues of matrix $g$. If
$|\lambda| =1, \lambda \ne 1$
(i.e. $g$ is an {\it elliptic} element of
$SU(1,1)$)
then $g$ has a fixed point on $D$ and hence
$\ell(g) =0$.
If
$\lambda =1$
(i.e. $g$ is {\it parabolic}) then $g$ have fixed points on circle
$|z| =1$
and
$\ell(g) =0$.

Let $g$ be {\it hyperbolic} (i.e.
$\lambda$
is real,
$\lambda \ne 1$).
Then there exist unique $g$-invariant geodesic
$L_g$
and the restriction of transformation
$z \mapsto z^{[g]}$
(see $\S$ \, 1) onto
$L_g$
be the shift of the line
$L_g$
on value
$|\ln(\lambda)|$.
Hence in this case

$$\ell(g) =| \ln(\lambda)|$$
Let
${{\cal M}}_n$
be a sphere with $n$ handles
$(n > 1)$.
Let
$\Gamma_n$
be fundamental group of ${{\cal M}}_n$. Recall that {\it Teichm\"uller space}
$Te_n$
is the space of all embeddings
$\Gamma_n \rightarrow \PSL_2(\R)$ such that image is a lattice.
Let
$\psi \in Te_n$.
Define {\it length function}
$\ell_{[\psi]}(q)$
on the group by the expression

$$\ell_{[\psi]}(q) = \ell(\psi(q))$$
This function has a simple geometric interpretation. Fix the embedding
$\psi : \Gamma_g \rightarrow \PSL_2(\R)$
and consider Riemann surfaces
$D/\psi(\Gamma_n)$.
Let
$h \in \Gamma_n$,
let
$L_{\psi(h)}$
be geodesic on $D$ invariant with respect to
$\psi(h)$.
Then the projection of
$L_{\psi(h)}$
onto
$D/\psi(\Gamma_n)$
is closed geodesic on
$D/\psi(\Gamma_n)$
and
$\ell_{[\psi]}(h)$
is the length
of this closed geodesics.

Hence we get that numbers
$\ell_{[\psi]}(h)$
(where
$h \in \Gamma_n)$
are lengths of closed geodesics on
$D/\psi(\Gamma_n)$.
It is well-known that the collection of lengths of closed geodesics determines the Riemann surface and hence the function
$\ell_{[\psi]}$
determines the point
$\psi \in Te_n$.

{\bf 3.6. Proof of Proposition 3.9.}
Let
$h \in \Gamma_n$.
Let
$\psi \in Te_n$.
Let $L$ be
$\psi(h)$
-invariant geodesic. Let
$u \in L$
be the point of $L$ nearest to $0$.

Evidently

$$\ell_{[\psi]}(h^n) = |n| \cdot \ell_{[\psi]} (h)$$
Now let us calculate asymptotics of
$\parallel \gamma (\psi(h^n))\parallel^2$
if
$n \rightarrow \infty$.
Using (1.4) we establish

$$\parallel \gamma (\psi(h^n))\parallel^2 = d(0,0^{[\psi(h^n)]}) + O(1)$$
Evidently

$$ |d(0,0^{[\psi(h^n)]}) - |n| \cdot \ell_{[\psi]} (h)|
= |d(0,0^{[\psi(h^n)]}) - d(u,u^{[\psi(h^n)]} | \le $$
$$ \le d(0,u) + d(0^{[\psi(h^n)]}, u^{[\psi(h^n)]}) = 2d(0,u) $$
Hence the expression

$$\parallel \gamma (\psi(h^n)) \parallel^2 -n \cdot \ell_{[i]} (h)$$
is bounded (for a given
$h \in \Gamma_n)$.
Bearing in mind that the length function
$l_{[\psi]}$
determines the embedding $\psi$ 
we apply lemma 2.1 and  finish the proof of the Proposition 3.9 in the case
$s_1 = s_2$.
It is easy to show that two lengths functions
$\ell_{[\psi]}(g), \ell_{[\psi^{\prime}]} (g)$
can't be proportional. This argument finishes the proof in the case
$s \ne s^{\prime}$.

\section{ Actions  of groups on $\R$-trees} 
\stepcounter{yyy}
\setcounter{xxx}{1}
\def\Vert{{\mbox{\rm Vert}}}
\def\Edge{{\mbox{\rm Edge}}}
\def\I{{\cal I}}

\def\risuno{
\unitlength0.8mm
\begin{picture}(80,100)
\put(15,0){
\put(0,0){\line(0,1){80}}
\put(4,0){\line(0,1){80}}
\put(10,0){\line(0,1){36}}
\put(8,46){\line(0,1){34}}
\put(72,0){\line(0,1){36}}
\put(76,0){\line(0,1){80}}
\put(82,0){\line(0,1){80}}
\put(70,46){\line(0,1){34}}
\put(10,36){\line(1,0){62}}
\put(4,42){\line(1,0){45}}
\put(76,40){\line(-1,0){45}}
\put(8,46){\line(1,0){62}}
\put(12,34){{\oval(8,8)[br]}
\put(3,-8){$A^+$}
\put(4,0){\vector(0,1){1}}}

 \put(12,48)
{{\oval(8,8)[tr]}
\put(3,8){$A^+$}
\put(0,4){\vector(-1,0){1}}}

 \put(-2,40)
{{\oval(8,8)[l]}
\put(-11,0){$A^+$}
\put(0,-4){\vector(1,0){1}}}
\put(25,42){
\put(0,0){\vector(0,1){4}}
\put(0,4){\vector(0,-1){4}}
\put(0,4){\linethickness{0.05mm}\line(0,1){10}}
\put(0,16){$a_+(v_\alpha,e_\gamma)$}
}
\put(27,42){
\put(0,0){\vector(0,-1){6}}
\put(0,-6){\vector(0,1){6}}
\put(0,-6){\linethickness{0.05mm}\line(0,-1){15}}
\put(0,-23){$a_-(v_\alpha,e_\gamma)$}
}
\put(52,40){
\put(0,0){\vector(0,1){6}}
\put(0,6){\vector(0,-1){6}}
\put(0,6){\linethickness{0.05mm}\line(0,1){20}}
\put(-15,28){$a_-(v_\beta,e_\gamma)$}
}
\put(55,40){
\put(0,0){\vector(0,-1){4}}
\put(0,-4){\vector(0,1){4}}
\put(0,-4){\linethickness{0.05mm}\line(0,-1){25}}
\put(-15,-34){$a_+(v_\beta,e_\gamma)$}
}
\put(4,75){
\put(0,0){\vector(1,0){4}}
\put(4,0){\vector(-1,0){4}}
\put(4,0){\linethickness{0.05mm}\line(1,0){10}}
\put(12,2){$a_-(v_\alpha,A^+(e_\gamma))$}
}
\put(40,30){$\bf e_\gamma$}
\put(-20,70){$\bf A^+_{v_\alpha}(e_\gamma)$}
\put(-20,10){$\bf A^-_{v_\alpha}(e_\gamma)$}
\put(76,40){\circle*{2}}
\put(84,40){$\bf v_\beta$}
\put(4,36){$\bf v_\alpha$}
\put(4,42){\circle*{2}}

}\end{picture}}

\def\risunokk{\unitlength0.8mm
\begin{picture}(80,100)
\put(15,0){
\put(0,0){\line(0,1){80}}
\put(4,0){\line(0,1){80}}
\put(10,0){\line(0,1){36}}
\put(8,46){\line(0,1){34}}
\put(72,0){\line(0,1){36}}
\put(76,0){\line(0,1){80}}
\put(82,0){\line(0,1){80}}
\put(70,46){\line(0,1){34}}
\put(10,36){\line(1,0){62}}
\put(4,42){\line(1,0){45}}
\put(76,40){\line(-1,0){45}}
\put(8,46){\line(1,0){62}}
\put(2,0){\line(0,1){80}}
\put(78,0){\line(0,1){80}}
\put(80,0){\line(0,1){80}}
\put(6,0){\line(0,1){40}}
\put(8,0){\line(0,1){38}}
\put(6,44){\line(0,1){36}}
\put(72,0){\line(0,1){36}}
\put(72,44){\line(0,1){36}}
\put(74,0){\line(0,1){38}}
\put(74,42){\line(0,1){38}}
\put(8,38){\line(1,0){66}}
\put(6,40){\line(1,0){62}}
\put(8,38){\line(1,0){66}}
\put(6,40){\line(1,0){66}}
\put(4,42){\line(1,0){70}}
\put(6,44){\line(1,0){66}}

\put(15,42){
\put(0,0){\linethickness{0.5mm}\vector(0,1){4}}
\put(0,4){\linethickness{0.05mm}\line(0,1){3}}
\put(0,9){$\Delta_+(v_\alpha,e_\gamma)$}
}
\put(20,42){
\put(0,-6){\linethickness{0.5mm}\vector(0,1){6}}
\put(0,-6){\linethickness{0.05mm}\line(0,-1){9}}
\put(0,-13){$\Delta_-(v_\alpha,e_\gamma)$}
}
\put(52,40){
\put(0,6){\linethickness{0.5mm}\vector(0,-1){6}}
\put(0,6){\linethickness{0.05mm}\line(0,1){3}}
\put(-5,13){$\Delta_-(v_\beta,e_\gamma)$}
}
\put(55,40){
\put(0,0){\linethickness{0.5mm}\vector(0,-1){4}}
\put(0,-4){\linethickness{0.05mm}\line(0,-1){25}}
\put(-15,-34){$\Delta_+(v_\beta,e_\gamma)$}
}
\put(4,75){
\put(4,0){\linethickness{0.5mm}\vector(-1,0){4}}
\put(4,0){\linethickness{0.05mm}\line(1,0){10}}
\put(12,2){$\Delta_-(v_\alpha,A^+(e_\gamma))$}
}
\put(4,62){
\put(0,0){\linethickness{0.5mm}\vector(-1,0){4}}
\put(0,0){\linethickness{0.05mm}\line(1,0){8}}
\put(7,2){$\Delta_+(v_\alpha,A^+(e_\gamma))$}
}
\put(4,30){
\put(0,0){\linethickness{0.5mm}\vector(1,0){6}}
}
\put(4,27){
\put(-4,0){\linethickness{0.5mm}\vector(1,0){4}}
}


}\end{picture}}

\def\risunokkk{\unitlength0.8mm
\begin{picture}(80,100)
\put(15,0){
\put(0,0){\line(0,1){80}}
\put(4,0){\line(0,1){80}}
\put(10,0){\line(0,1){36}}
\put(8,46){\line(0,1){34}}
\put(72,0){\line(0,1){36}}
\put(76,0){\line(0,1){80}}
\put(82,0){\line(0,1){80}}
\put(70,46){\line(0,1){34}}
\put(10,36){\line(1,0){62}}
\put(4,42){\line(1,0){45}}
\put(76,40){\line(-1,0){45}}
\put(8,46){\line(1,0){62}}
\put(2,0){\line(0,1){80}}
\put(78,0){\line(0,1){80}}
\put(80,0){\line(0,1){80}}
\put(6,0){\line(0,1){40}}
\put(8,0){\line(0,1){38}}
\put(6,44){\line(0,1){36}}
\put(72,0){\line(0,1){36}}
\put(72,44){\line(0,1){36}}
\put(74,0){\line(0,1){38}}
\put(74,42){\line(0,1){38}}
\put(8,38){\line(1,0){66}}
\put(6,40){\line(1,0){62}}
\put(8,38){\line(1,0){66}}
\put(6,40){\line(1,0){66}}
\put(4,42){\line(1,0){70}}
\put(6,44){\line(1,0){66}}

\put(76,15){
\put(0,0){\circle*{1}}
\put(6,0){\circle*{1}}
\put(17,8){$x=z_0$}
\put(16,8){\linethickness{0.05mm}\vector(-1,0){10}}
\put(8,0){$u_0$}
\put(0,0){\linethickness{0.05mm}\line(-1,0){6}}
\put(-9,2){$v_1$}
\put(0,0){\linethickness{0.3mm}\line(1,0){6}}
}
\put(20,40){
\put(0,0){\circle*{1}}
\put(0,2){\circle*{1}}
\put(0,0){\linethickness{0.3mm}\line(0,1){2}}
\put(0,0){\linethickness{0.05mm}\line(0,-1){6}}
\put(0,2){\linethickness{0.05mm}\line(0,1){6}}
\put(1,8){$v_2$}
\put(1,-8){$u_1$}
}

\put(2,70){
\put(0,0){\circle*{1}}
\put(2,0){\circle*{1}}
\put(0,0){\linethickness{0.3mm}\line(1,0){2}}
\put(0,0){\linethickness{0.05mm}\line(-1,0){5}}
\put(0,0){\linethickness{0.05mm}\line(1,0){7}}
\put(8,1){$u_2$}
\put(-8,2){$v_3$}
\put(-10,-10){\linethickness{0.05mm}\vector(1,0){8}}
\put(-17,-8){$y=z_3$}

}

\put(40,56){
\put(-1,-1){$z_2$}
\put(4,0){\linethickness{0.05mm}\vector(1,0){30}}
\put(0,-4){\linethickness{0.05mm}\vector(0,-1){10}}
\put(-4,0){\linethickness{0.05mm}\vector(-1,0){32}}
}

\put(56,10){
\put(-1,-1){$z_1$}
\put(4,0){\linethickness{0.05mm}\vector(1,0){16}}
\put(0,4){\linethickness{0.05mm}\vector(0,1){26}}
\put(-4,0){\linethickness{0.05mm}\vector(-1,0){46}}
}

}\end{picture}}

{\bf 4.1 $\R$-trees,}
(see [CM], [M], [MS], [Par], [Sha1], [Sha2]).
A {\it metric tree} is a tree ( a tree is a graph without cycles) with given lengths of edges. We will consider a metric tree as 
a metric space. Denote by $\Vert(I)$ the set of all vertices
of the metric tree $I$. Denote by $\Edge(I)$ the set of
all edges of the metric tree $I$.
 Consider a chain

$${\cal T}_{(1)} \stackrel{\sigma_1}\longrightarrow {\cal T}_{(2)} \stackrel{\sigma_2}  \longrightarrow {\cal T}_{(3)} \stackrel{\sigma_3}  \longrightarrow \ldots \no\label
{ 4.1}$$
where
${\cal T}_{(i)}$
are metric trees and
$\sigma_i$
are isometric embeddings. We mention that the image of
a vertex of a tree ${\cal T}_{(j)}$ can be not a vertex of the
tree ${\cal T}_{(j+1)}$! 
 Consider the inductive limit $\I$ of this chain. Metric
spaces which can be obtain by such a way are called
{\it $\R$-trees}.

 We say that a point of $I$ is
a {\it vertex} if it is a vertex of some tree $\I_{(j)}$.
The set of vertices of $\R$-tree $\I$ we denote by
$\Vert (I)$.

  Generally speaking {\it $\R$-trees have no edges.}

{\bf Remark 4.1.} Our definition of $\R$-trees slightly differs
from standard definition. Our $\R$-tree is official
$\R$-tree with countable number of vertices.

Consider a $\R$-tree $\cal T$.  For each two points
$x,y \in {\cal T}$
there exist a unique isometric embedding of some interval
$[0,\alpha] \subset \R$
to ${\cal T}$ such that
$\lambda (0) =x, \lambda (\alpha) = y$
(of course
$\alpha$
equals to distance
$d(x,y)$
between $x$ and $y$). The image
$\lambda ([0,\alpha])$ of the segment $([0,\alpha])$ 
is denoted by
$[x,y]$.

 Let
$\Gamma$
be a group. An {\it 
action} of
$\Gamma$
{\it on
$\R$
-tree $\I$} is (by definition) an action by isometric transformation of
$\R$
-tree $\I$.

{\bf 4.2. Examples of actions of groups on metric trees. }

{\bf Example 4.2.}
Consider a lattice
$\Gamma \subset \PSL_2(\K)$
where $\K$ is a locally compact nonarchimedian field. The group
$\PSL_2(\K)$
acts on Bruhat-Tits tree
${\cal T}_q$
(see subsection 1.6) and hence
$\Gamma$
also acts on
${\cal T}_q$.

{\bf Example 4.3. Trees associated to free groups.}
Remind the standard construction of tree associated to free
group.
Denote by
$F_n$
the free group with $n$ generators
$a_1, \ldots, a_n$.
 The vertices of the tree
${\cal T}(F_n)$
are the elements of the group
$F_n$.
Two vertices
$w, w^{\prime}$
of
${\cal T}(F_n)$
are incident to the same edge if

$$w^{\prime} = w\cdot a_j^{\pm 1}$$
The left action of the group
$F_n$
on itself bornes the action of
$F_n$
on
${\cal T}(F_n)$. Of course this tree be the usual Cayley graph of the group $F_n$.
Let us equip the graph
${\cal T}(F_n)$
by some additional structures

The first, we have $n$ different types of edges corresponding to the generators
$a_1,\ldots, a_n$.
We will call them by
$\{a_j\}$
-type edges. The second, let us fix the canonical orientation on each edge. Consider an edge $(w,wa_j) $
where $w$ is an irreducible word and the last letter
of $w$ is not $a_j^{-1}$. Then by definition $w$ is the origin 
of the edge and
$wa_j$
is the end.

Now we can say that the length of
$\{a_j\}$
-type edge equals to
$s_j > 0$.
Then we obtain the metric tree which we denote by
${\cal T}(F_n; s_1, \ldots s_n)$.

{\bf Example 4.4. Action of free groups on Bruhat-Tits tree
 ${\cal T}_\infty$.} Consider the preceding example.
Let us contract all edges of the types
$\{a_2\}, \ldots, \{a_n\}$
(or equivalently let us consider the graph
${\cal T}(F_n;1,0,\ldots, 0)$.
We obtain the new oriented graph
${\cal T}^{*}(F_n)$.
The vertices of
${\cal T}^{*}(F_n)$
are enumerated by cosets
$F_n/F_{n-1}$
where
$F_{n-1}$
is the free group generated by
$a_2,\ldots, a_n$.
Edges
$[w,w^{\prime}]$
of
${\cal T}^{(*)}(F_n)$
have the form

$$ w = \psi_1 a^{\varepsilon_1}_1 \ldots \psi_{k}^{\varepsilon_{k}}a_1^{\varepsilon_k} F_{n-1} $$

$$w^{\prime} = \psi_1 a^{\varepsilon_1}_1 \ldots \psi^{\varepsilon_{k}}_{k} a^{\varepsilon_{k}}_1 \varphi \, a_1 F_{n-1} \no\label {f1} $$
where

$$\varepsilon_i = \pm 1,\qquad \psi_i, \varphi \in F_{n-1}$$
and the word $w$ is irreducible ($w$ is the origin,
$w^{\prime}$ is the end). Each vertex of
${\cal T}^{*}(F_n)$
is the origin of countable number of edges and the end of countable number of edges. The action of the group
$F_n$
on $\Edge({\cal T}^*(F_n)$
is transitive and the stabilizer of each edge is trivial. The bijection

$$\mbox{\rm Edge} \, ({\cal T}^{*} (F_n)) \leftrightarrow F_n \no\label {f2}$$
is given by the formula

$$[w,w^{\prime}] \leftrightarrow \psi_1a^{\varepsilon_1}_1 \ldots \psi_{k} a^{\varepsilon_{k}}_{k} \varphi$$
where
$[w,w^{\prime}]$
is the edge (4.2).

{\bf Example 4.5. Actions of surface groups,} see [MSh].
Let $M_n$ be a sphere with $n$ handles ($n>1$). Let $\Gamma_n$
be its fundamental group. Consider a finite family of closed
non-self-intersecting and not pairwise intersected curves 
$C_1,\dots, C_l$ on $M_n$.
Assume $M_n\setminus \{ C_1,C_2,\dots,C_l\}$ is a union
of spheres with holes. Consider the universal covering $D$
of $M_n$. Remind that $D$ is the Lobachevskii plane.
Then we obtain a countable number of infinite curves $S_1,S_2\dots$
on $D$ which cover curves $ C_1,C_2,\dots,C_l$. These curves separate Lobachevskii plane $D$ to countable number of pieces $U_1,U_2,\dots$.
 
Now we will construct a tree. Edges of this tree are
numerated by curves $S_1, S_2,\dots$ and vertices are
numerated by open domains $U_1,U_2,\dots$. A
 vertex $U_\alpha$ and a edge $S_\beta$ are
 incindent iff the curve $S_\beta$  lies on the boundary of
the domain $U_\alpha$ .

{\bf 4.3. Examples of actions of groups on $\R$-trees.}

{\bf Example 4.6. Train-tracked graphs}, see [Fock],
on train tracks see[PH].
Consider a graph $\J$. Denote by $\Edge^*(\J)$ the set
of all pairs
$$(v_\alpha,e_\gamma)\in \Vert(\J)\times\\Edge(\J)$$
such that vertex $v_\alpha$ is incident to the
edge $e_\gamma$. In other words $\Edge^*(J)$ is the set of
 oriented edges. A {\it train-tracked graph} is a graph
$\J$ equipped with the following structures:

  1. For each vertex $v_\alpha$ we fix a cyclic order
on the set of edges incident to $v_\alpha$. We denote by
$A^+_{v_\alpha}(e_\gamma)$ the edge which follows to
$e_\gamma$ with respect to this cyclic order. By
$A^-_{v_\alpha}(e_\gamma)$ 
we denote the preceding edge.

  2. We fix two function $a_+(v_\alpha,e_\gamma)$ 
and $a_-(v_\alpha,e_\gamma)$ on the set
$\Edge^*(J)$ with values in nonnegative numbers. These
functions satisfies two conditions.

$$\qquad{\mbox{\rm a)}} \qquad\qquad\qquad \qquad a_+(v_\alpha,e_\gamma)=
 a_-(v_\alpha,A^+_{v_\alpha}(e_\gamma))$$

b) Let a edge $e_\gamma$ is incindent to vertices
$v_\alpha$ and $v_\beta$. Then

$$a_+(v_\alpha,e_\gamma)+a_-(v_\alpha,e_\gamma)=
a_+(v_\beta,e_\gamma)+a_-(v_\beta,e_\gamma)\no\label{Train}$$

{\bf Remark 4.7}.
It is natural to imagine that edges of our graph are strips
and the expression (4.4) is the width
of the strip. It is also natural to think that
strip from the edge $e_\alpha$ is continued to the adjacent edges $A^+_{v_\alpha}(e_\gamma)$ and $A^-_{v_\alpha}(e_\gamma)$ and separated between them, see picture 1.

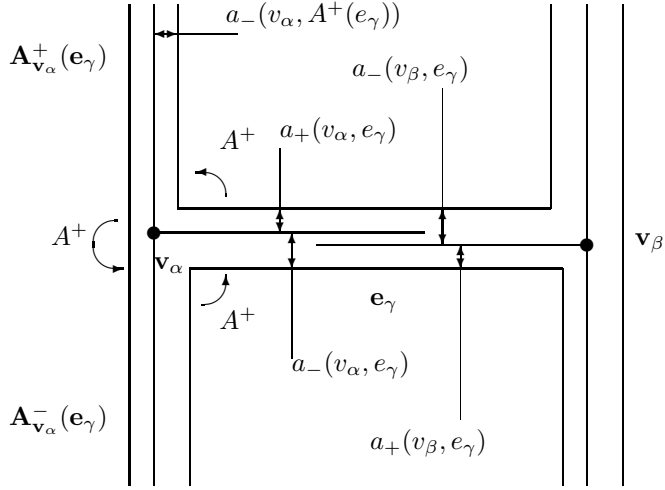
\risuno

{\it The universal covering} of a train-tracked graph is
defined by the obvious way. We obtain a tree with
train-tracked structure.

Now for each a train-tracked tree $\J$ we will construct
some $\R$-tree ${{\cal R}}(\J)$. For each element
 $(v_\alpha,e_\gamma)\in\Edge^*(J)$
we consider two directed segments
$\Delta_+(v_\alpha,e_\gamma)$ and
$\Delta_-(v_\alpha,e_\gamma)$ 
having length
$a_+(v_\alpha,e_\gamma)$ and
$a_-(v_\alpha,e_\gamma)$. Consider the disjoint union
$\Omega$ of all
oriented segments $\Delta^\pm(v_\alpha,e_\gamma)$.

Now for each pair $(v_\alpha,e_\gamma)\in\Edge^*(J)$
we will identify  segments $\Delta^+(v_\alpha,e_\gamma)$ and $\Delta^-(v_\alpha,A^+_{v_\alpha}(e_\gamma))$
by isometrical orientation inversion way (remind that
their lengths are equal). Consider arbitrary edge
$e_\alpha$ let $v_\alpha,v_\beta$ be its end. Then we identify
the following segments 

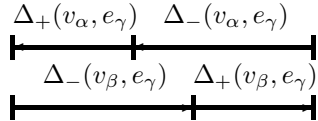
\begin{picture}(100,40){\unitlength0.8mm\linethickness{0.5mm}
\put(20,2){
\put(0,0){\vector(1,0){30}}
\put(30,0){\vector(1,0){20}}
\put(20,10){\vector(-1,0){20}}
\put(50,10){\vector(-1,0){30}}

\put(0,-2)
{\put(0,0){\line(0,1){4}}
\put(30,0){\line(0,1){4}}
\put(50,0){\line(0,1){4}}
\put(0,10){\line(0,1){4}}
\put(20,10){\line(0,1){4}}
\put(50,10){\line(0,1){4}}}

\put(0,4){
\put(0,10){$\Delta_+(v_\alpha, e_\gamma)$}
\put(25,10){$\Delta_-(v_\alpha, e_\gamma)$}
\put(5,0){$\Delta_-(v_\beta, e_\gamma)$}
\put(30,0){$\Delta_+(v_\beta, e_\gamma)$}
}

}}\end{picture}

(with inversion of direction).

 After all these gluings we obtained some
topological space ${{\cal R}}(\J)$, see picture 2
(elements of the quotient a branching thicklines
on the picture).

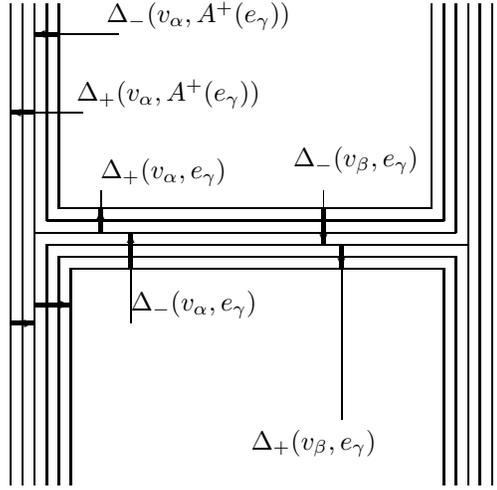
\risunokk

Now we will define a metric
$d(x,y)$ on the quotient space ${{\cal R}}(\J)$. Consider
all finite chains
$$x=z_0,z_1,z_2,\dots,z_m=y\no\label{train}$$
such that for each $j$ there exist a segment
$\Delta^\pm(v_\alpha,e_\gamma)$ such that elements
$z_j$ and $z_{j+1}$ of quotient have representatives
$u_j$ and $v_{j+1}$ in the segment $\Delta^\pm(v_\alpha,e_\gamma)$.

Then
$$d(x,y)=\min (\sum_j  |u_j-v_{j+1}|)$$
there minimum is given by all chains (4.5),
see picture 3.

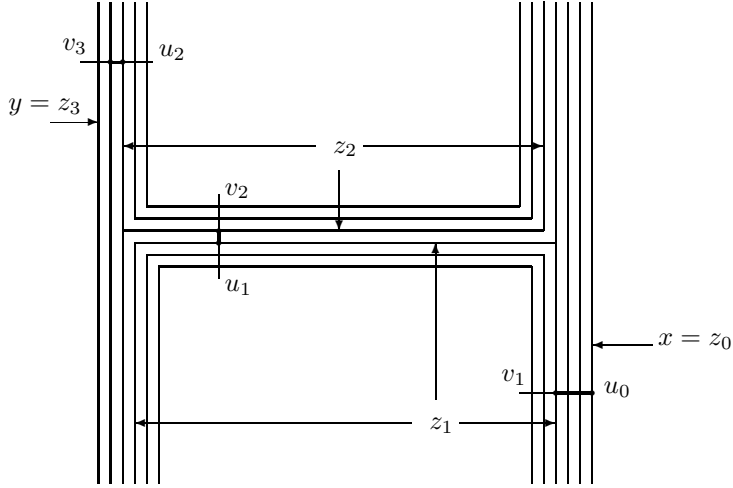
\risunokkk

It is easy to see that the space ${{\cal R}}(\J)$ is a $\R$-tree.

\def\Lat{{\rm Lat}}
{\bf Example 4.7. Bruhat-Tits trees over nonlocal fields},
see [Tit].
Let $K$ be a valuated nonarchimedian field. Let $\Lambda
\subset\R$ be the group of valuation (remind the definition
of {\it valuation}: It is a function
 $\lambda:K\setminus 0\rightarrow\Lambda$ such that 
$\lambda (xy)=\lambda(x)+\lambda(y)$ and
$\lambda (x+y)\le \{\max (\lambda(x),\lambda(y)\}$).
Consider the space $\Lat$ of 
all lattices in $K^2$. For each lattice $R$ we
define function $\gamma_R:K^2\rightarrow \Lambda$
by the condition:
$\gamma(v)$ is the minimum of $\lambda(x)$ by all
$x\in K$ such that
$x^{-1}v\in R$.

 Let $\Lat/K^*$ be space of lattices
defined up to dilatation (see ss.1.5)  Define distance $d(R_1,R_2)$ 
 between points
$R_1,R_2\in\Lat/K^*$  by the following rule:
$$d(R_1,R_2)=\max_{v\in\K^2}(\gamma_{R_1}(v)-\gamma_{R_2}(v))-
\min_{v\in\K^2}(\gamma_{R_1}(v)-\gamma_{R_2}(v))$$.

Points of  $\Lat/K^*$ are vertices of our $\R$-tree.
(we omit obvious construction of tree itself).

Obviously group $GL_2(K)$ acts on our tree.

{\bf 4.5. Length function.}

{\bf Proposition 4.8}
(see [CS]).{\it  Let
${\cal T}$
be an 
$\R$-tree and
$q : {\cal T} \rightarrow {\cal T}$
be an isometric map. Then there exist two possibilities:

 1) There exists a fixed point
$x \in {\cal T} (qx =x)$.

 2) There exists a $q$-invariant line
$C_q$
and the restriction of $q$ to
$C_q$
is a shift of
$C_q$.}

The line
$C_q$
is called the {\it axis} of isometry $q$.

{\bf Theorem 4.9}
(see [CM], [MS]) {\it Let a group 
$\Gamma$
acts on an
$\R$
-tree
${\cal T}$.

 a) If
$\Gamma$
is finitely generated then there exists a minimal 
$\Gamma$-
invariant subtree
$\I \subset {\cal T}$

 b) If there exists
$g \in \Gamma$
such that the map
$g : {\cal T} \rightarrow {\cal T}$
has no fixed points then there exists a minimal
$\Gamma$
-invariant subtree $\I$. }

If
$\I = {\cal T}$
then the action is called {\it minimal}.

Some conditions on action of group in this theorem are really
nessesary (countrexample: action of group of $p$-adic 
matrices having a form 
$\left(\begin{array}{cc}1& \ast\\0&1\end{array}\right)$
on Bruhat--Tits tree).

 Consider an action of a group
$\Gamma$
on an
$\R$
-tree $\cal T$. {\it  Length function} (or {\it translation length function})
$\ell : \Gamma \rightarrow \R$
is defined by

$$\ell (g) = \min_{x \in {\cal T}}
 (x,gx)$$
In fact a minimum exists (see Proposition 4.8). If $g$ has a fixed point, then
$\ell (g) =0$.
If $g$ has no fixed point then for each
$x \in C_g$
we have
$\ell (g) = d(x,gx)$.

{\bf Theorem 4.10.}{\it
(see [CM]) Consider minimal actions of the group
$\Gamma$
on trees
${\cal T}_1, {\cal T}_2$
with the same length function
$\ell$.
Let
$\ell$
be not identical zero on the commutant
$[\Gamma, \Gamma]$.
Then these two actions are isomorphic.}

(If a group $\Gamma$ have a fixsd point on the absolut
of tree then $\ell$ equals to 0 on $[\Gamma,\Gamma]$.)

Denote by
$\P  (\Gamma)$
the space of real functions on the discrete group
$\Gamma$
defined up to multiplier. Let us equip the space
$\P  (\Gamma)$
by the topology of pointwise convergence.

There exists a simple list of properties defining length functions, see [Chi],[CM], [Par]. This list allows to prove
the following theorem

{\bf Theorem 4.11}{\it
(see [CM]) The space of all length functions on
$\Gamma$
is compact in
$\P  (\Gamma)$.}

{\bf 4.6. Thurston compactification.}
Consider the Teichm\"uller space
$Te_n$,
(see subsection 3.5), let
$\psi \in Te_n$ be a embedding of the surface group
$\Gamma_n$ to $\PSL_2(\R)$
and let
$\ell_{[\psi]}$
be its length function (see subsection 3.5). Let us consider this function as a point of the space
$\P  (\Gamma_n)$.
Hence we obtained the embedding
$Te_n \rightarrow \P (\Gamma_n))$.
Denote by
$\overline{Te_n}$
the closure of
$Te_n$
in
$\P (\Gamma_n)$.
This closure is called the {\it Thurston compactification} of the Teichm\"uller space (see [Thu], [MS]).

{\bf Theorem 4.12.}{\it
(see [MS]) Each points of the boundary
$\overline{Te_n}\setminus Te_n$
is a length function for some isometric action of the group
$\Gamma_n$ on a $\R$-tree.}

{\bf Remark 4.13.} Moreover the space of homomorphisms
$\Gamma_n \rightarrow \PSL_2(\R)$
(= the space of actions of
$\Gamma_n$
on Lobachevskii plane) has the natural compactification by actions of
$\Gamma_n$
on
$\R$
-trees. The same is true for the space of homomorphisms from a given finitely generated discrete group
$\Gamma$
to
$\SO(1,n)$.
It is valid also for homomorphisms
$\Gamma$
to
$\PSL_2(\K)$
(the latter is corollary of the theorem 4.11).

{\bf 4.7. Affine hilbert space connected with $\R$-tree.}
 We will construct by arbitrary
$\R$
-tree
${\cal T}$
some Hilbert space
$H({\cal T})$.

To each pair
$x,y$
of vertices of
${\cal T}$
we associate the formal vector
$e(x,y)$.
We assume
$e(x,y) = -e(y,x)$.
Let
$x,y,u,v$
be vertices of
${\cal T}$.
Assume

$$\langle e(x,y), e(u,v)\rangle = \pm \quad \mbox{\rm (length of intersection} \quad [x,y] \cap [u,v]) \no\label{ 4.2}$$
We choose the sign ``+'' if segments
$[x,y]$
and
$[u,v]$
are directed to the same side on the intersection
$[x,y] \cap [u,v]$.
We choose ``-'' in the opposite case.

Consider the Hilbert space
$H({\cal T})$
such that

 a)
$e(x,y)$
are elements of
$H({\cal T})$
and their scalar products are given by the formula 
(4.6). 

 b) Linear combinations of vectors
$e[x,y]$
are dense in
$H({\cal T})$.

{\bf Remark 4.14.}
Let
$x,y,z$
be vertices of
${\cal T}$.
Then

$$e(x,y) + e(y,z) + e(z,x) =0 \no\label{ 4.3}$$
Indeed if
$z \in [x,y]$
(or
$x \in [y,z]$
or
$y \in [x,z]$)
it follows from (4.6). Otherwise there exists unique vertex
$\underline{u}$
such that
$u \in [x,y] \cap [y,z] \cap [x,z]$:

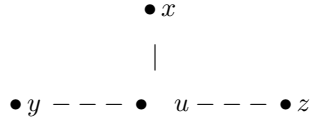
$$\bullet \, x $$
$$\vert \, \, $$
$$\bullet \, y \, --- \, \bullet \quad u --- \, \bullet \, z$$
Then

$$e(x,y) = e(x,u) + e(u,y) $$
$$e(y,z) = e(y,u) + e(u,z) $$
$$e(z,x) = e(z,u) + e(u,x) $$
Adding this equalities and bearing in mind equality
$e(\alpha,\beta) = -e(\beta,\alpha)$
we obtain (4.7).

{\it Uniqueness} of
$H({\cal T})$
is obvious. Let us prove the {\it existence }  of
$H({\cal T})$.
Consider  a metric tree
${\cal J}$.
Fix a vertex $a\in \Vert{\cal J}$. This vertex defines a orientation on the tree by the following rule:

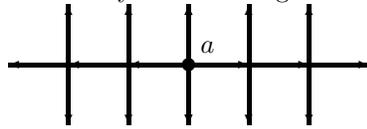
\begin{picture}(60,20){\linethickness{0.5mm}
\put(30,10){
\put(0,0){\circle*{2}}
\put(2,2){$a$}
\multiput(0,0)(10,0){3}
{\put(0,0){\vector(1,0){10}}
\put(0,0){\vector(0,1){10}}
\put(0,0){\vector(0,-1){10}}  }

\multiput(0,0)(-10,0){3}
{\put(0,0){\vector(-1,0){10}}
\put(-0,0){\vector(0,1){10}}
\put(-0,0){\vector(0,-1){10}}  }}

}\end{picture}

Formally a edge $[x,y]$ is positive directed
if $x\in[a,y]$. For each positive directed edge
 $[x_j,y_j]$ we associate a formal vector $e[x_j,y_j]$.
We assume  $e[x_j,y_j]$
is a orthogonal basis and
 $$\|e[x_j,y_j]\|^2=\mbox{\rm length of} [x_j,y_j]$$.
We define the space $H({\cal T})$ as the space generated by vectors $e[x_j,y_j]$.

Now consider a chain

$${\cal T}_{(1)} \stackrel{\sigma_1}  \longrightarrow {\cal T}_{(2)} \stackrel{\sigma_2}  \longrightarrow {\cal T}_{(3)} \stackrel{\sigma_3}  \longrightarrow \ldots \no$$

For this chain   we associate the chain

$$H({\cal T}_1) \rightarrow H({\cal T}_2) \rightarrow \ldots \no\label{ 4.4}$$

The space
$H({\cal T})$
is an inductive limit of chain (4.9) of hilbert spaces(i.e. completion of
$\cup_k H({\cal T}_k)$).

Let a group
$\Gamma$
acts on the
$\R$
-tree
${\cal T}$.
This action induces the unitary representation of
$\Gamma$
in
$H({\cal T})$
by the rule

$$\pi (g)e(x,y) = e(gx,gy) \no\label{ 4.5}$$

{\bf Lemma 4.15.}
{\it Fix a vertex
$x_0 \in \Vert({\cal T})$.
Then the formula

$$\Pi_s(g)v = \pi(g)v + s\cdot e(x_0,g x_0) \no $$
defines the isometric action of
$\Gamma$
in
$H({\cal T})$.}

{\bf Proof.}
We have to check the equality

$$e(x_0,g_1g_2x_0) = \pi (g_1)e(x_0,g_2x_0) + e(x_0,g_1x_0)$$
Right side equals to

$$e(g_1x_0, g_1g_2x_0) + e(x_0,g_1 x_0)$$
and we apply (4.7).

{\bf Remark 4.16.} If to consider another origin $x_0\in \Vert({\cal T})$
we obtain an equivalent affine isometric action.

{\bf Remark 4.17.}
Let
$\theta$
be a point of the Thurston compactification
$\bar Te_m$ of Teichmuller space $Te_m$.
We canonically associate to
$\theta$
an AI-irreducible affine isometric action
$\Pi^{\theta} = (\pi_{\theta},\gamma_{\theta})$
of
$\Gamma_m$. In the
the case
$\theta \in Te_m$
the action
$\Pi^{\theta}$
is (by definition) the restriction of real affine action (see subsection 2.5) of
$\PSL_2(\R)$
to
$\Gamma_m$.
If
$\theta \in \bar Te_m \setminus Te_m$
then
$\Pi_{\theta}$
is an AI-irreducible subrepresentation of the action 4.6.

There arises the natural question: Is the map
$\theta \mapsto \Pi^{\theta}$
continuous? I don't know answer on this question, but most likely the answer is affirmative. Observe that in both cases
$(\theta \in Te_m$
or
$\theta \in \bar Te_m\setminus Te_m)$
the length function
$\ell_{\theta}$
and the conditionally positive defined function
$\parallel \gamma_{\theta}(g)\parallel^2$
(see below $\S$ \, 5)
are related by equality

$$\ell_{\theta}(g) = {\lim}_{n \rightarrow \infty}\frac 1 n \parallel \gamma (g^n) \parallel^2 \no $$

{\bf 4.8. Existence of AI-irreducible subrepresentation.}

{\bf Theorem 4.18.}{\it
Consider a minimal action of a group
$\Gamma$
on an
$\R$-tree
${\cal T}$
satisfying the conditions of theorem 4.10.
Then there exist
AI-irreducible subrepresentations in
$H({\cal T})$.}

{\bf Proof.}
We have to prove that the linear representation
$\pi(\cdot)$
of the group
$\Gamma$
in
$H({\cal T})$
given by the formula (4.10) doesn't contains weakly the trivial representation of
$\Gamma$ (see proposition 2.9).
Assume that it is not so. Consider a sequence
$h_j \in H({\cal T})$
such that
$\langle h_j,h_j\rangle =1$
and
$\langle \pi(g) h_j,h_j \rangle \rightarrow 1$
for each
$g \in G$.
Consider an element
$g \in G$
which has no fixed point on
${\cal T}$ (i.e. $\ell(g)\neq 0$.
Let
$C_g$
be the axis of $g$. Let

$$h_j = p_j +q_j, p_j \in H(C_g), q_j \in H({\cal T} \setminus C_g).$$
Then

$$\langle \pi(g)h_j,h_j\rangle = \langle \pi(g) (p_j + q_j), p_j + q_j)\rangle = \langle \pi(g)p_j,p_j \rangle$$
(since
$H(C_g)$
is orthogonal to
$H({\cal T}\setminus C_g)$
and
$H({\cal T}\setminus C_g)$
is the orthogonal sum of spaces
$H(Q_{\alpha})$
where
$Q_{\alpha}$
are connected components of the graph
${\cal T}\setminus C_q)$. We have

$$\langle \pi (g)h_j,h_j\rangle \le \parallel p_j\parallel^2 =1 - \parallel q_j\parallel^2$$
Hence $\|q_j\|\rightarrow 0$. In other words
$\rho (h_j,H(C_g)) \rightarrow 0$
for
$j \rightarrow \infty$, where $\rho$ denote distance between point and subspace.
Hence we obtain that for each finite collection
$g_1, \ldots, g_{\varkappa} \in G$

$$\rho(h_j, 
\cap_{\alpha =1}^{\varkappa} H (C_{g_{\alpha}})) \rightarrow 0 \no\label{exis}$$
if
$j \rightarrow \infty$.
Consider elements
$g_{\alpha}$
such that
$I=\cap_{\alpha =1}^{\varkappa} H (C_{g_{\alpha}})) $
is a finite segment (otherwise $\Gamma$ has fixed points on the boundary of tree). For big $N$ we have
$g^N_1I \cap I =\varnothing$
and hence
$g^N_1(\cap_{\alpha =1}^{\varkappa} H (C_{g_{\alpha}}))$
is orthogonal to
$\cap_{\alpha =1}^{\varkappa} 
H(C_{g_{\alpha}})$.
This and
$(4.13)$
contradicts to the condition
$\langle \rho (g_1)h_j,h_j\rangle \rightarrow 1$.

{\bf Proposition 4.19.}{\it
Consider two nonequivalent actions of  group
$\Gamma$
satisfying to the conditions of theorem 4.10. Then corresponding AI-irreducible subrepresentations of
$\Gamma$
are not equivalent.}

{\bf Proof.}
It is sufficient to reconstruct the length function of
the action
$\Pi$
(see Theorem 4.10). Denote
$e(x_0,gx_0)$
by
$\gamma (g)$.
Then

$$\parallel \gamma (g)\parallel^2 = d(x_0,gx_0).$$
Let $u$ be the element of the axis
$C_g$
nearest to
$x_0$.
Then

$$d(x_0,g^n x_0) = d(x_0,u) +d(u,g^nu) + d(g^nu,g^nx_0) = 2d(x_0,u) + n \ell (g)$$

$$u \qquad gu \qquad \qquad \qquad \qquad \qquad \qquad \qquad \qquad \qquad \qquad \qquad \qquad g^nu
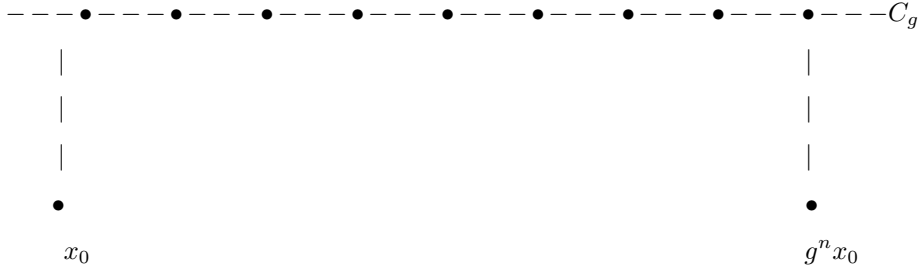
$$
$$--- \bullet --- \bullet --- \bullet --- \bullet --- \bullet --- \bullet --- \bullet --- \bullet --- \bullet --- C_g$$
$$\mid \qquad \qquad \qquad \qquad \qquad \qquad \qquad \qquad \qquad \qquad \qquad \qquad \qquad \qquad \mid \qquad$$
$$\mid \qquad \qquad \qquad \qquad \qquad \qquad \qquad \qquad \qquad \qquad \qquad \qquad \qquad \qquad \mid \qquad$$
$$\mid \qquad \qquad \qquad \qquad \qquad \qquad \qquad \qquad \qquad \qquad \qquad \qquad \qquad \qquad \mid \qquad$$
$$\bullet \qquad \qquad \qquad \qquad \qquad \qquad \qquad \qquad \qquad \qquad \qquad \qquad \qquad \qquad \bullet \qquad$$
$$\qquad x_0 \qquad \qquad \qquad \qquad \qquad \qquad \qquad \qquad \qquad \qquad \qquad \qquad \qquad \quad g^nx_0 \qquad$$
Hence (for the given $g$)

$$\parallel \gamma (g^n)\parallel^2 -n \ell (g) = O(1)$$
Hence the function
$\parallel \gamma (g)\parallel^2$
determines the length function.$\blacksquare$

{\bf Example 4.20.} Consider the Bruhat-Tits tree
${\cal T}_n$. Then our hilbert space is the space 
$l_2^\pm(\Edge^*({\cal T}_n))$ in notations of ss.1.5.
Obviously action of $\Aut({\cal T}_n)$ in 
$l_2^\pm(\Edge^*({\cal T}_n))$ 
is AI-reducible. Its AI-irreducible component
can be realized in the space $\ker(\blacktriangledown)=
Im({{\cal P}})$.

{\bf 4.9. Irreducibility theorem.}
Let
${\cal T}$
be a metric tree . Let a group
$\Gamma$
acts on
${\cal T}$.
For a vertex
$x_0 \in {\cal T}$
denote by
$\Gamma (x_0)$
the stabilizer of
$x_0$
in
$\Gamma$.
For an edge
$\delta$
denote by
$\Gamma[\sigma]$
the stabilizer of
$\sigma$
in
$\Gamma$.

Recall that for a metric tree
${\cal T}$
the Hilbert space
$H({\cal T})$
is identified with the space of functions on the set 
$\Edge({\cal T})$
provided by the scalar product

$$\langle f_1,f_2\rangle = \sum_{\sigma_j \in \mbox{\rm Edge} ({\cal T})} f_1(\sigma_j) \overline{f_2(\sigma_j)} \cdot \mbox{\rm length} \, (\sigma_j)$$

{\bf Theorem 4.21.}{\it
Consider an action of the group
$\Gamma$
on an oriented metric tree
${\cal T}$.
Let
$\Gamma$
be transitive on the set of vertices. Let
$x_0$
be a vertex of
${\cal T}$.
Assume that for each edge
$\sigma$
incindent to
$x_0$
the index of
$\Gamma[\sigma]$
in
$\Gamma (x_0)$
is infinite. Then the canonical affine isometric representation of
$\Gamma$
in
$H({\cal T})$
is AI
-irreducible.}

{\bf Proof.}
Suppose that there exist an invariant subspace $K$ in
$H({\cal T})$.
Let
$Af = Tf + \tau$
be the projector to $K$ (where $T$ is a linear operator in
$H({\cal T})$
and
$\tau = \tau(\sigma) \in H({\cal T})$
is a function on Edge
$({\cal T})$).
The operator $A$ commutes with the affine isometric action

$$\Pi(g) f(\sigma) = f(g\sigma) + \gamma_g(\sigma)$$
where
$\sigma \in \mbox{\rm Edge}({\cal T}), \gamma_g(\cdot) \in H({\cal T})$.
Hence

$$T\, \gamma_g(\sigma)-\gamma_g(\sigma) = \tau (g\sigma) - \tau (\sigma) \no\label{irre}$$
Let
$g \in \Gamma (x_0)$.
Then the left part of (4.4) equals 0 (since
$\gamma_g(\sigma) =0$ far all $\sigma\in\Edge$).
Hence the function
$\tau (\sigma)$
is
$\Gamma (x_0)$
-invariant. But all orbits of
$\Gamma (x_0)$
on Edge
$({\cal T})$
are infinite. Bearing in mind that
$\tau (\cdot) \in l_2$
we obtain 
$\tau =0$.
Hence
$T \, \gamma_g(\sigma) = \gamma_g(\sigma)$.
Bearing in mind that
$\Gamma$ 
is transitive on the set of vertices we obtain that a linear span of vectors
$\gamma_g$
is dense in
$H({\cal T})$
and hence
$T = E$.$\blacksquare$

Consider the action of free group $F_n$ on 
the Bruhat-Tits tree ${\cal T}_\infty$ constructed
in Example 4.4.

{\bf Corollary 4.22.}{\it
The canonical affine isometric action of
$F_n$
in
$H({\cal T}^{*} (F_n))$
is irreducible.}

{\bf Proof.}
Indeed this action satisfies the conditions of Theorem 4.21. 

The bijection (4.3) identifies the space
$H({\cal T}^{(1)}(F_n))$
and
$l_2(F_n)$.
Rewrite the canonical action of
$F_n$
in
$H({\cal T}^{(1)} (F_n))$
in terms of the space of
$l_2(F_n)$.
Consider

$$g = \psi_1 a^{\varepsilon_1}_1 \psi_2 a^{\varepsilon_2} \ldots \psi_{k}a^{\varepsilon_{k}} \psi_{k +1} \in F_n \no$$
where
$\psi_i \in F_{n-1}, \varepsilon_i = \pm 1$
and the word (.) is irreducible. Then

$$\Pi (g) f(h) = f(gh) + \gamma_g (h) \no$$
where
$h \in F_n$
and the function
$\gamma_g(h)$
is given by the formula

$$\gamma_g(h) = \left\{\begin{array}{c}\varepsilon_j \quad \mbox{\rm if} \quad h = \psi_1a^{\varepsilon_1}_1 \ldots \psi_j a^{\varepsilon_j} \\
0 \quad \mbox{\rm otherwise} \end{array} \right.$$

{\bf Remark 4.23.}
Recall that the regular representation of
$F_n$
is a factor-representation, nevertheless its affine ``deformation'' is AI-irreducible.

{\bf Example 4.24.}
Consider once more the tree
${\cal T}(F_n)$ (see Example 4.3).
Let us contract all edges of types
$\{a_{\alpha +1}\}, \ldots, \{a_n\}$.
Denote the tree which we obtained by
${\cal T}^{(\alpha)}(F_n)$.
The vertices of
${\cal T}^{(\alpha)} (F_n)$
are enumerated by cosets
$F_n/F_{n-\alpha}$
and the group
$F_n$
has
$\alpha$
orbits on the set of edges (the stabilizer of each edge is trivial). If
$\alpha < n$
the conditions of the theorem 4.21 are fulfiled and hence the affine isometric action of
$F_n$
in the space
$H({\cal T}^{\alpha}(F_n))$
is AI
-irreducible. Observe also that
$H({\cal T}^{(\alpha)}(F_n))$
can be identified by the natural way with
$\alpha$
copies of the space
$l_2(F_n)$.

The affine isometric action of
$F_n$
in
$H({\cal T}(F_n,s_1, \ldots, s_n))$
is reducible (see the arguments of Example 4.20).

{\bf 4.10. Bibliographical remarks.} Bruhat--Tits trees
appeared in [BH] (1967) as a partial case of Bruhat--Tits
buildings (ensembles). Serre used Bruhat--Tits trees 
as tool for investigation of discrete groups ([Ser1-2]). $\R$-Trees
were discovered in [Chi] and [Tit] (1976--1977).

  Construction of affine isometric actions and
conditionally positive definite functions associated to
metric tree is old and I don't know when it appeared,
see [Ols2], [Mar],III.3. Spherical functions (they are closely
related to affine isometric actions) on Bruhat--Tits tree of
nonarchimedian nonlocal field were classified 
in [Ism1]-[Ism2]. Conditionally positive definite
functions associated to $\R$-tree were mentioned in [HV].
On $\R$-trees see reviews [Sha1-2].

\section{General remarks (continuation)} 
\stepcounter{yyy}
\setcounter{xxx}{1}
This section is continuation of preliminaries.

{\bf 5.1. Boson Fock space}, see [Ber1],[Ner2].
{\it Boson Fock space}
$F_n$
{\it with $n$ degrees of freedom} is the Hilbert space of holomorphic functions on
$\C^n$
with scalar product

$$\{ f_1,f_2\} = \int_{\C^n} f_1(z) \overline{f_2(z)} \exp (-|z|^2) \prod \frac{dz_i \overline{dz_i}} {\pi}$$
The function
$f(z) =1$
is called by {\it vacuum vector}, and its norm equals 1. We have the canonical isometric embedding
$I_n : F_n \rightarrow F_{n+1}$,
if
$f(z_1, \ldots, z_n) \in F_n$
then

$$(I_nf)(z_1,\ldots, z_{n+1}) = f(z_1, \ldots , z_n) \in F_{n+1}$$
Consider the chain of embeddings

$$F_0 \rightarrow F_1 \rightarrow F_2 \rightarrow \ldots$$
Its inductive limit
$F_{\infty}$
in the category of Hilbert spaces (i.e. the completion of
$\cup F_j$)
is called by {\it Boson Fock space with infinite number degrees of freedom}. Elements of
$F_{\infty}$
can be considered as holomorphic functions on
$l_2$, see [Ber1], [Ner2].

Let $H$ be a complex separable Hilbert space. Then $H$ is isomorphic to
$l_2$.
Hence elements of
$F_{\infty}$
can be considered as holomorphic functions on $H$. 
We denote this space of holomorphic functions  by
$F(H)$.

{\bf 5.2 Fock representation of a group of isometries.}
Denote by Isom
$(H)$
the group of all affine isometric transformations of a Hilbert space $H$.

Consider an element
$\sigma \mapsto A\sigma + \gamma$
of the group Isom
$(H)$.
Consider the operator

$$\Exp(A,\gamma)f(z) = f(Az + \gamma) \exp (\langle z,A^{-1}\gamma\rangle - \frac 1 2 \langle \gamma,\gamma\rangle) \no\label {5.1}$$
in Boson Fock space
$F(H)$.
The following theorem is well-known (and simple).

{\bf Theorem 5.1.}{\it
$\Exp(A,\gamma)$
is a unitary projective representation of the group Isom
$(H)$.
Moreover

$$\Exp(T,\gamma)\Exp(T^{\prime},\gamma^{\prime}) = \exp \{i \, \Im < \gamma, (T^{\prime})^{-1} \gamma'\} \cdot
\Exp(T^{\prime\prime}, \gamma^{\prime\prime}) \no\label 
{5.2}$$
(where
$T^{\prime\prime} = TT^{\prime}, \gamma^{\prime\prime} = T^{\prime} \gamma + \gamma^{\prime})$.}

{\bf 5.3 Araki scheme.}
Let $G$ be a group and let
$\Pi_s = (\pi,s\cdot \gamma)$
be its affine isomorphic action. Then the formula

$$\rho_s = \Exp \circ \Pi_s$$
defines the unitary (projective) representation of the group $G$ in the boson Fock space
$F(H)$.

{\bf Remark 5.2}
Let $K$ be a {\it real} Hilbert space and let
$\Pi = (\pi,\gamma)$
be an affine isometric action of $G$ in $K$. Consider complexification
$K_{\C}$
and representation
$\Exp$
of Isom
$(K_{\C})$.
Then (see (5.2))
$Exp \circ \Pi$
is linear (nonprojective !) representation of $G$.

Application of this construction to locally compact groups gives extremely reducible representations of $G$ (the boson Fock space is so large and the group $G$ is so small!). In the next subsection we briefly discuss some interpretations of this construction.

For infinite dimensional groups such constructions are more interesting, see [Ara], [VGG1-2], [Ism3-5], [Ols2], [Ols3], [Ner 1-2].

{\bf 5.4 Conditionally positive definite functions,} 
(see [Scho], [Gui1], [Mar], [PS] [HV]). Let $G$ be a topological group. Denote by
${\cal P}(G)$
the set  which elements are triples
$(H,\pi,\sigma)$
where

 i) $H$ is a Hilbert space

 ii) $\pi$
is a unitary representation of $G$ in $H$

 iii) $\sigma \in H$
is $G$-cyclic vector.

We say that triples
$(H,\pi,\sigma)$
and
$(H^{\prime},\pi^{\prime}, \sigma^{\prime})$
are  equivalent if there exists a $G$-intertwining unitary operator
$U : H \rightarrow H^{\prime}$
such that
$U \, \sigma = \sigma $. Let
$P = (H,\pi,\sigma), P^{\prime}=(H^{\prime},\pi^{\prime},\sigma^{\prime}) \in {\cal P}(G)$. Let us define the product
$PP^{\prime} = (H^{\prime\prime}, \pi^{\prime\prime}, \sigma^{\prime\prime})$.
For this consider the triple

$$(H \otimes H^{\prime}, \pi \otimes \pi^{\prime}, \sigma \otimes \sigma^{\prime})$$
Let
$H^{\prime\prime}$
be $G$-cyclic span of the vector
$\sigma \otimes \sigma^{\prime}, \pi^{\prime\prime}$
be restriction of representation
$\pi \otimes \pi^{\prime}$
on
$H^{\prime\prime}$
and
$\sigma^{\prime\prime} = \sigma \otimes \sigma^{\prime} \in H^{\prime\prime}$.
It is well-known that the semigroup
${\cal P}(G)$
is isomorphic to the (multiplicative) semigroup
${\cal R}(G)$
of all positive definite functions on $G$ (Gelfand--Naimark--Segal construction, see [Dix],13.4).

Recall that a continuous function
$\psi : G \rightarrow \C$
is called {\it positive definite} if for each
$g_1, \ldots, g_n \in G$
the hermitian form

$$Q(h_1, \ldots, h_n) = \sum_{i,j =1}^n \psi (g^{-1}_ig_j)h_ih_j$$
is nonnegative defined (see [Dix], 13.4). Let
$(H,\pi,\sigma) \in {\cal P}(G)$.
Then

$$\psi(g) = \langle \pi(g)\sigma, \sigma\rangle_H$$
is an element of
${\cal R}(G)$
and the correspondense
${\cal P}(G) \leftrightarrow {\cal R}(G)$
is bijective.

We are interested by 1-parametric semigroups in
${\cal P}(G)$, i.e. families
$P(t) \in {\cal P}(G), t\ge 0$
such that

$$P(t_1+t_2) = P(t_1) P(t_2)$$

{\bf Example 5.3.}
Let
$G = \widetilde{\SL_2(\R)}$
be the universal covering of
$\SL_2(\R)$.
Let
$\pi_t$
be the representation of $G$ in the space of holomorphic functions in the disc
$|z| < 1$
given by the formula

$$\pi_t \left(\begin{array}{cc} a & b \\ \bar b & \bar a \end{array}\right) f(z) = f\left( \frac {az+b} {\bar bz + \bar a}\right) (\bar b z + \bar a)^{-t}$$
If
$t > 0$ then
the representation
$\pi_t$
is unitary in some scalar product. Let
$H_t$
be the corresponding Hilbert space and
$\sigma_t=\sigma_t(z) =1$.
Then

$${\cal P}(t) = (H_t,\pi_t,\sigma_t)$$
be the one-parametric subsemigroup in
${\cal P}(\overline{\SL_2(\R)})$.
We also have

$$\langle \pi_t \left(\begin{array}{cc} a & b \\ \bar b & \bar a \end{array}\right) 1, 1\rangle = \bar a^{-t}
$$

Consider an one-parametric semigroup
$P(t)$
in
${\cal P}(G)$.
Let
$\psi_t$
be the corresponding family
${\cal R}(G)$.
Evidently
$\psi_t(g)$
has a form

$$\psi_t(g) = \exp (-t\cdot \varphi(g))$$
where
$\varphi$
is a function on $G$. There arises the natural question:
describe functions
$\varphi(g)$ such that
functions 
$\exp(-t\cdot \varphi (g))$
are positive definite for all
$t > 0$ ? The answer is well-known (see [Scho]),
 this are the so-called conditionally positive definite functions. Recall that a continuous function
$\varphi$
on $G$ is called {\it conditionally positive definite (or negative definite)} if for each
$g_1,\ldots , g_n \in G$
the hermitian form

$$Q(h_1, \ldots, h_n) = \sum_{i,j}\phi(g^{-1}_i g_j) h_ih_j$$
is nonpositive definite on the hyperplane

$$h_1 + \ldots + h_n =0.$$

{\bf Remark 5.4.}
 The description of all conditionally positive definite functions on
$\R$
is given by the classical Levy-Khinchin theorem (see [Shi],III.5). For nonabel groups the question on 
analogy of Levy-Khinchin theorem have two completely
different versions: the question on infinitely divisible elements
in semigroup ${{\cal M}}(G)$ of probabilistic measures
(see [Hey]) and question on infinitely divisible representations.

{\bf 5.5. Conditionally positive definite functions and affine isometric actions.}
Consider an affine isometric action of a group $G$ on a 
{\it real} Hilbert space $H$:

$$\Pi(g)v = \pi(g) v + \gamma (g) \no\label{ 5.3}$$
It is easy to prove that

$$\varphi (g) = \parallel \gamma (g) \parallel^2$$
is a conditionally positive definite function satisfying the conditions

 a) $\varphi(g)$
is real 

 b) $\varphi (e) =0$

 c) $\varphi (g) = \varphi (g^{-1})$

{\it Inversely}, let
$\varphi$
be a conditionally positive definite function satisfying the conditions a) - c). Then where exist a unique {\it real} action (5.3) such that

 a) $\varphi (g) = \parallel \gamma (g) \parallel^2$

 b) There are no $G$-invariant closed {\it linear} subspaces in $H$.

The construction is simple. Conditions (0.3)-(0.4)
imply
$$
\|\gamma(g^{-1})\|=\|\gamma(g)\|^2+\|\gamma(h)\|^2
-2\langle \gamma(h),\gamma(g)\rangle$$
and hence we have to construct a hilbert space $H$ and a system of vectors $\gamma(h)$ such that
$$\langle \gamma(h),\gamma(g)\rangle=
\frac 12 (\phi(g^{-1}h)-\phi(g)-\phi(h))$$

We omit details, see for instance [HV].

It is easy to describe a one-parametric semigroup
$P(t) \in {\cal P}(G)$
corresponding to the semigroup
$\exp(-t \parallel \gamma (g) \parallel^2)$
in
${\cal R}(G)$.
For this consider the representation
$\Exp \circ \Pi_{\sqrt t}$
(see subsection 5.3) and the cyclic span of vacuum vectors.

\section{  Fock representation  of the  semigroup of probabilistic measures on a group} 

\stepcounter{yyy}
\setcounter{xxx}{1}

{\bf 6.1 Construction of representation}
Let $G$ be a topological group. Denote by
${\cal M}_0(G)$
the semigroup of probabilistic measures on $G$ with compact support. Denote by
$\mu * \nu$
the convolution of measures. Consider an affine isometric 
action of $G$ in a Hilbert space $H$:

$$\Pi(g)v = \pi(g) v + \gamma (g)$$

{\bf Theorem 6.1}{\it
(see [Ner 2]). Let
$\mu \in {\cal M}_0(G)$.
Consider the operator
$L(\mu)$
in the boson Fock space
$F(H)$
defined by the formula

$$L(\mu)f(z) = \exp (d + \langle u,c\rangle )f(Az+b)$$
where

$$ A = A(\mu) = \int_G \pi (g) d\mu \no\label {6.1} $$
$$b = b(\mu) = \int_G \gamma(g) d\mu) \no\label{ 6.2} $$
$$c = c(\mu) = \int_G \gamma (g^{-1})d\mu \no\label {6.3}$$
$$d = d(\mu) = -\int_G \parallel \gamma (g) \parallel^2 d\mu \no\label {6.4}$$
Then

 a)

$$L(\mu)L(\nu) = \exp (i \cdot \Im \langle b(\nu), c(\mu)\rangle) \cdot L(\mu * \nu) \no\label {6.5}$$

 b)
$\parallel L(\mu) \parallel \le 1$
for each 
$\mu$.

Hence we get a projective representation of the semigroup
${\cal M}_0(G)$.
We also get central extensions of
${\cal M}_0(G)$
given by the  formula (6.5).}

{\bf 6.2 Examples of central extensions.}
Let
$\Pi = (\pi,\gamma)$
be a {\it complex} affine isometric action (see subsection 2.5) of a group $G$. Then (by Theorem 6.1) we obtain the central 
extension of semigroup
${\cal M}_0(G)$.
The elements of this extension are pairs
$(\mu,x) \in {\cal M}_0(G) \times \R$
and the product is given by the formula

$$(\mu,x)\circ (\nu ,y) = (\mu *\nu,x+y+\sigma (\mu, \nu))$$
where

$$ \sigma (\mu,\nu) = \Im \langle b(\nu), c(\mu) \rangle = 
 \Im \int_G \int_G \langle \gamma (g), \gamma(h^{-1})\rangle d\mu (g) d\nu (h) $$

{\bf Example 6.2.} Consider the case
$G = \PSL_2(\R)$.
The affine isometric action of $G$ is defined in
$\S$ \, 1. When using the formula (1.2) we obtain

$$\sigma(\mu,\nu) =-\Im \int_{\PSL_2(\R)} \int_{\PSL_2(\R)} \ln \left( \frac { a(g_1g_2)} {a(g_1)a(g_2)} \right) d\mu (g_1) d\nu (g_2)\no\label{ts1}$$

{\bf Example 6.3.}Let now
$G = \Gamma_n$
be the fundamental group of a sphere with $n$ handles. Consider an embedding
$\psi : \Gamma_n \rightarrow \PSL_2(\R)$.
Elements of
${\cal M}_0(\Gamma_n)$
are finite linear combinations
$\Sigma \, p_j \delta (g_j)$,
where
$p_j>0, \Sigma \, p_j =1$
and
$\delta (g_j)$
is the unit measure supported in the point 
$g_j \in \Gamma$.
Then the cocycle is given by the formula

$$\sigma(\sum^M_{j=1} p_j\delta(g_j), \sum^N_{k =1} 
q_{k}\delta (h_{k})) = \Im \, \sum\, \ln 
\frac{a(\psi(g_jh_{k}))} 
{a(\psi(g_j))a(\psi(h_{k}))} 
p_jq_{k}\no\label{ts2}$$

{\bf Example 6.4.} Let
$G \simeq\Z^{2n}$
be a lattice in
$\C^n$.
The cocycle on
${\cal M}_0(G)$
is given by the formula

$$ \sigma (\sum^N_{j=1} \alpha_j\delta (\sigma_j), \sum^M_{k =1} \alpha^{\prime}_{k} \delta (\sigma^{\prime}_{k})) 
= \sum^N_{j=1} \sum^M_{k =1} \alpha_j \alpha^{\prime}_{k} \Im \, \langle \sigma_j, \sigma^{\prime}_{k} \rangle $$

{\bf 6.3 Irreducibility.}

{\bf Theorem 6.2.} {\it
Let
$\Pi = (\pi,\gamma)$
be an irreducible isometric affine action of a group $G$. Let the unitary representation
$\pi$
doesn't contain weakly the trivial representation of $G$. Then the representation
$L(\mu)$
of the semigroup
${\cal M}_0(G)$
is irreducible.}

{\bf Proof.}
(I am using the argument communicated to me by Ismagilov). Denote by $H$ the space of representation 
$\Pi$.
Consider a measure
$\mu \in {\cal M}_0(G)$
such that
$\parallel \pi (\mu)\parallel < 1$, see proposition 2.7.
Denote by
$\mu^\circ$
the image of
$\mu$
with respect to a map
$g \mapsto g^{-1}$
from $G$ to itself. Let
$\nu = \mu^\circ * \mu$
(hence
$\nu^\circ = \nu)$.
Consider the map
$\tilde \Pi(\nu) : H \rightarrow H$ given by the formula

$$\tilde \Pi (\nu)v = \int_G \pi(g)v d\nu (g) + \int_G \gamma (g) d\nu(g)$$
Let
$w \in H$
be the fixed point of the map    $\tilde \Pi(\nu)$.
Without loss of generality we assume
$w =0$.
Hence

$$0 = \int_G \gamma(g) d\nu(g) = \int_G \gamma (g^{-1})d\nu(g^{-1}) = \int_G \gamma (g^{-1})d\nu (g)$$
and vectors
$b(\nu)$
and
$c(\nu)$
defined by (6.2)-(6.3) equals zero. Hence
$L(\nu)$
has the form

$$L(\nu) f(z) = \lambda \cdot f(A(\nu)z)$$
where
$\lambda \in \C$,
and
$(A(\nu))^* = A \, (\nu^\circ) = A(\nu)$.
Hence the vacuum vector 1 is the unique maximal vector of the operator
$L(\nu)$.
Hence each
${\cal M}_0(G)$
-intertwining operator fixes the vacuum vector. It remains to prove that vector 1 is
${\cal M}_0(G)$
-cyclic. It is quite simple, for details see [Ner 2],10.2.

\section{Groups of automorphisms of bundle} 

\stepcounter{yyy}
\setcounter{xxx}{1}
\def\Ams{{\rm Ams}}

{\bf 7.1.
The group}
$\B(G)$.
Let $X$ be a Lebesgue measure space with continuous measure
$\mu$
such that
$\mu (X) =1$.
Recall that all such spaces are isomorphic to the segment
$[0,1]$ equipped by the usual Lebesgue measure.
Denote by
$\Ams(X)$
the group of all measure preserving maps
$X \rightarrow X$.

Let $G$ be a topological group. Denote by
${{\cal F}}(X,G)$
the group of all functions
$h : X \rightarrow G$
such that there exist a compact set
$K  \subset G$ (subset $K$ depends on function $h$)
satisfying the condition
$\mu \{x \in X : h(x) \notin K\} =0$.
Denote by
$$\B(G) = \B(X,G)$$
the semidirect product of the group
$\Ams (X)$
and the group
${{\cal F}}(X,G)$.
Elements of
$\B(G)$
are pairs
$$(p(x),h(x)) \in \Aut (x) \times {{\cal F}}(X,G)$$
and the product is given by the formula

$$ (p_1(x),h_1(x))(p_2(x),h_2(x))  
=(p_1(p_2(x)); h_1(p_2(x)) \cdot f_2(h_2(x)) $$

{\bf 7.2. Induced actions of the group $\B(G)$.}
Let a group $G$ acts by transformations
$\varkappa (g)$
of a set $M$. Then the group
$\B(G)$
act on
$X \times M$
by transformations

$$\varkappa^\circ(p,h)(x,m) = \left(p(x), h(x)m\right)$$

{\bf 7.3. Induced representations of the group $\B(G)$.}
Consider a unitary  representation $\rho$
 of a group $G$ in a Hilbert space $H$. Consider the space
$L^2(X,H)$
of all
$L^2$-functions
$X \rightarrow H$,
the scalar product in
$L^2(X,H)$
is given by the formula

$$\langle f_1,f_2 \rangle = \int_X (f_1(x), f_2(x))_Hd\mu (x)$$
Define the unitary representation
$\rho^\circ$
of
$\B(G)$
in
$L^2(X,H)$
 by the formula

$$\rho^\circ(p,h) f(x) = \rho(h(x))f(p(x))$$
where
$p \in \Ams(X), h\in {{\cal F}}(X,G)$.

{\bf 7.4. Induced affine isometric actions of the group
$\B(G)$.}
Consider an affine isometric action

$$\Pi(g)v= \pi(g) v + \gamma (g)$$
of a group $G$ in a Hilbert space $H$. Define the affine isometric action
$\Pi^\circ = (\pi^\circ,\gamma^\circ)$
the group
$\B(G)$
in the space
$L^2(X,H)$
by the formula

$$\Pi^\circ(p,h)f (x) = \pi(h(x))f(p(x)) + \gamma (h(x)) \no\label {7.1}$$

{\bf 7.5. Fock representations of the group $\B(G).$
(``{\it Araki multiplicative integral}'')
}
Let a group $G$ has a nontrivial affine isometric action
$\Pi = (\pi,\gamma)$.
Consider the affine isometric action
$\Pi^\circ = (\pi^\circ,\gamma^\circ)$
given by the formula (7.1). Restriction of  the representation
$\Exp$
of the group
$\Isom(L^2(X,H))$
to the subgroup
$\B(G)$
give  projective unitary representation
$R = \Exp \circ \Pi^\circ$
of the group
$\B(G)$
in the Fock space 
$F(L^2(X,H))$
given by the formula

$$ R(p,q) f(w) = 
f(\Pi^\circ(p,q)w) \times$$ $$\times
 \exp \left[-\frac 1 2 \int_X \parallel \gamma (q(x)) \parallel ^2 d\mu (x) 
-\int_X \langle \pi (q(x))w(p(x)), \gamma (q(x)) \rangle_H d\mu(x))
\right]$$
where
$w = w(x) \in L^2 (X,H), f= f(w) \in F(L^2(X,H))$.
Using (5.2) we obtain

$$R(p,q) R(\tilde p, \tilde q) = \sigma ((p,q), (\tilde p, \tilde q)) \cdot
 R((p,q) \circ (\tilde p, \tilde q)) $$
where

$$ \sigma ((p,q), (\tilde p, \tilde q)) = \exp \left\{i \, \Im \int_X \langle \pi (\tilde q(x)) \cdot
 \gamma (q (\tilde p(x)), \gamma (\tilde q(x))\rangle_H d\mu(x)\right\} \no\label {7.2}$$
{\it Araki multiplicative integral} is the restriction
of this representation to the subgroup ${{\cal F}}(X,G)$.
On Araki multiplicative integral  see [VGG1-2], [Ber2], [Del1-2].

{\bf 7.6 Central extensions of $\B(G)$.}
Let a group $G$ has a nontrivial action
$\Pi = (\pi,\gamma)$
on a {\it complex} Hilbert space $H$. Then in the preceding subsection we constructed the central extension
$\widetilde{\B(G)}_\Pi$.
The elements of this extension are triples
$(p,q,a) \in \Ams (X) \times {\cal F} (X,G) \times \R$
and the multiplication is given by the formula

$$(p,q,a) \circ (\tilde p, \tilde q, \tilde a) = \left(\tilde p \circ p, (\tilde q \circ p) \cdot \tilde q, a + \tilde a +  \lambda ((p,q);(\tilde p, \tilde q)\right)$$
where
$\lambda (\ldots)$
is the expression in curly brackets in (7.2).

{\bf Example 7.1.}
Let
$G = \PSL_2(\R)$.
Then with the notation  of subsection 1.2 we have

$$\lambda (p,q); (\tilde p, \tilde q)) =-\Im \int_X \ln \left( \frac {a(\tilde q(p(x))q(x))} {a(\tilde q(p(x))a(q(x))} \right) d\mu (x)$$

Let
$\Gamma$
be a lattice in
$\PSL_2(\R)$.
Then the central extension of
$\B(\Gamma)$
is given by the same formula

{\bf Remark 7.2.}
The constructions of $\S$ \, 6 and $\S$ \, 7 are the parts of a more general construction obtained in [Ner 1], see also [Ner 2] (representations of the category of $G$-stochastic kernels).

{\bf 7.7. Group of automorphisms of $G$-bundle.} Let $G$ be a group.
Let $X=M$ be a compact smooth manifold equipped by some
smooth volume form $\omega$ (of course $M$ is a Lebesgue measure space).
 Consider a vector
$G$-bundle $B$ over a manifold $M$. Denote by
$\pi$
the projection of the total space $B$ to $M$. Let
$r : B \rightarrow B$
be a measurable transformation. We say that $r$ is an {\it measurable automorphism of the bundle} if

 a) $r$ maps fibres to fibres. Hence $r$ induces the map
$r^{\square} : M \rightarrow M$

 b) the map
$r^{\square}$
preserves the measure on $M$
 
 c) for each
$m \in M$
the map of fibres

$$r : \pi^{-1}(m) \mapsto \pi^{-1}(r^{\square}(m))$$
 is a element of the group 
$G$.

It is easy to see that the group of all {\it measurable} automorphisms of the bundle $B$ is isomorphic to the group
$\B(M,G)$.

Denote by $\B^\infty[B]$ the group of {\it smooth} automorphisms
of the fibre bundle $B\rightarrow M$.

Obviously
$$\B^\infty[B]\subset\B(M,G)$$.
In many cases $\B^\infty[B]$ is dense in the group $\B(M,G)$ in some
natural sense and this implies that the restriction of
Fock representation of  $\B(M,G)$ to $\B^\infty[B]$ is irreducible.

Now we obtain a possibility to embed different groups $\frak G$ to the group
$\B^\infty[B]$ and by such a way to construct unitary representations and central extensions of $\B^\infty[B]$ .

\section{ Constructions of representations 
 of groups of diffeomorphisms of two-dimensional manifold}
\stepcounter{yyy}
\setcounter{xxx}{1} 
\def\SDiff{{\rm SDiff}}

Here we want to construct several series of representations 
of groups of volume preserving diffeomorphisms
of two-dimensional manifolds.

{\bf 8.1
Construction with Jacoby matrix},see [Ner1]. Let
$M^2$
be a compact two-dimensional manifold, let $\omega$ be a volume form on
$M^2$. Let
$G = \SDiff \, (M^2)$
be a group of all volume preserving and orientation
preserving diffeomorphisms of
$M^2$.
Consider a volume preserving and orientation preserving smooth maps from some bounded domain
$\Omega \subset \R^2$
to some open dense domain on
$M^2$.
For each 
$q \in \SDiff \, (M^2)$
consider a corresponding piecewise smooth map
$\tilde q : \Omega \rightarrow \Omega$.
Denote by
$D(\tilde q(x))$
the Jacoby matrix of
$\tilde q(x)$
in the point 
$x \in \Omega$.
Then the formula

$$q \mapsto (\tilde q; D(\tilde q(x)) \no$$
defines the embedding
$\SDiff (M^2) \rightarrow \B(\PSL_2(\R))$.
The image of
$\SDiff(M^2)$
is dense in
$\B(PSL_2(\R))$
(see [Ner 1]).

Restricting the Fock representation (see subsection 7.5) of
$\B(PSL_2(\R))$
to
$\SDiff (M^2)$
we obtain irreducible representation of
$\SDiff (M^2)$.
Analogically the central extension of
$\B(\PSL_2(\R))$
(see subsection 7.6) induces the nontrivial central extension of
$\SDiff(M^2)$.

{\bf Remark 8.1}. Consider the tangent
$\PSL_2(\R)$
-bundle
$B \rightarrow M^2$.
The group
$\SDiff (M^2)$
acts in $B$ by the obvious way. The formula (8.1) is ``coordinatization'' of this action.

{\bf Remark 8.2}
This construction  of representation of $\SDiff(M^2)$ 
can be repeated for the group of volume preserving diffeomorphisms of 2-dimensional $p$-adic manifolds.
 I don't know has it some sense or not.

{\bf 8.2 Ismagilov's embedding,}
see [Ism 2]. Let $M$ be a connected and  non simply connected compact manifold. Let $\omega$ be a volume form on $M$. Denote by
$\SDiff(M)$
the group of all diffeomorphism of $M$ preserving $\omega$. Let $\D(M)$ be the connected component of
$\SDiff(M^2)$
containing unity. Assume (for simplicity) that the fundamental group
$\Gamma$ of $M$ has a trivial center.

Let
$\tilde M$
be the universal covering over $M$. Then
$\pi : \tilde M \rightarrow M$
is a principle
$\Gamma$
-bundle over $M$. Let
$q \in \D(M)$.
Then there exists a unique covering diffeomorphism
$\tilde q : \tilde M \rightarrow \tilde M$
$(\pi \circ \tilde q = q \circ \pi)$.
Hence we get a (continuous) action of the group 
$\D(M)$
in $\Gamma$-principle bundle
$\pi : \tilde M \rightarrow M$.

So we obtain the embedding

$$\D(M) \rightarrow \B(M,\Gamma)\no\label{ism}$$

{\bf Theorem 8.3.}{\it
([Ism 2]) The image of
$\D(M)$
in
$\B(M,\Gamma)$
is dense.}

Now we will give another description of the same embedding. For this consider a compact fundamental domain
$M^* \subset \tilde M$
(i.e.
$\bigcup_{g_i \in \Gamma} g_i M^* = \tilde M; M^* \cap g_i M^*$
has zero measure for all
$g_i \ne e$).
Let
$q \in \D(M)$.
Define the function
$\lambda_q (x) : M^* \rightarrow \Gamma$
by the condition

$$\lambda_q(x) = h \Leftrightarrow \tilde q x \in hM^*$$
Then the embedding
$\D(M) \rightarrow \B(M^*, \Gamma)$
is given by the formula

$$q \mapsto (\tilde q, \lambda_q(x)).$$
By remarks from ss.7.2-7.3 each action of
$\Gamma$
induces the action of
$\D(M)$,
each unitary representation of
$\Gamma$
induces the unitary representation of
$\D(M)$.
Each affine isometric action of
$\Gamma$
induces an affine isometric action of
$\D(M)$
and representation of the group
$\D(M)$
in the boson Fock space.

{\bf 8.3. Action of group $\D$ on 3-dimensional manifolds.}
Let
$M_{n}$
be a sphere with 
$n$
handles
$(n > 1)$.
Then its universal covering
$\tilde M_{n}$
is the Lobachevskii plane
$\H^2$.

 Consider an embedding
$\lambda : \Gamma \rightarrow \PSL_2(\R)$.
The group 
$\PSL_2(\R)$
acts on the circle
$S^1$.
Let us restrict this action to the subgroup
$\Gamma$.
Concider embedding (8.2) and
 the induced action (see ss.7.2) of the group
$\D(M_{n})$.
In other words consider the
$S^1$-bundle
$\phi : R \rightarrow M_n$
induced from the principal
$\Gamma$
-bundle
$\H \rightarrow M_{n}$
and consider the action of
$\D(M_{n})$ 
in the space $R$. The space $R$ is a 3-dimensional manifold, orbits of the group
$\D(M)$
in $R$ are 2-dimensional submanifolds of $R$ isomorphic to
$\H^2$.
Each $\D(M)$-orbit
$\cal O$
is dense in $R$ (and the intersection of
$\cal O$
with each fibre
$\phi^{-1}(m) \simeq S^1$
is a dense subset in the fibre).

Consider the natural unitary representation of
$\D(M)$
in
$L^2(R)$:

$$L_\lambda(q)f(x) = f(qx) \cdot I(q,x)^{1/2}$$
where
$x \in R, q \in \D(M), I(q,x)$
is the Jacobian of the map
$x \mapsto qx$.

{\bf Proposition 8.4.}
{\it The representation
$L_{\lambda}$
is irreducible. If the embeddings
$\lambda, \lambda^{\prime} : \Gamma \rightarrow \PSL_2(\R)$
are not conjugate then the representations
$L_{\lambda}, L_{\lambda^{\prime}}$
are not equivalent.}

{\bf Proof.}
It is the partial case of the Proposition 8.7 (see below). This corresponds to one of the points of principal series representation of
$\PSL_2(\R)$.$\blacksquare$

Let $\pi$ be an irreducible unitary representation of
$\PSL_2(\R)$.
Consider the representation 
$L_{\lambda,\pi}$
of
$\D(M)$
induced from the restriction  $\pi{\LARGE |}\Gamma$.

{\bf Proposition 8.5.}{\it

 a) Let $\pi$  be not a representation of discrete series. Then
$L_{\lambda,\pi}$ 
is irreducible.

 b) Let
$\pi,\pi^{\prime}$
are not representation of discrete series. Then
$L_{\lambda,\pi} = L_{\lambda^{\prime},\pi^{\prime}}$,
only in the case
$\pi = \pi^{\prime}$
and
$\lambda$
is conjugate to
$\lambda^{\prime}$.

 c) Let $\pi$ be a representation of discrete series. Then
$L_{\lambda,\pi}$
be a type ${\bf II}_1$ factor-representation.}

{\bf Proof.} It is a corollary of theorem 8.5 and 
theorems 3.1-3.3.
$\blacksquare$

{\bf 8.4.Induced affine actions.}
Let us consider an affine action
$\Pi_s$
of
$\PSL_2(\R)$
(see $\S$ \, 1) and its restriction to
$\Gamma$.
Then we obtain a representation of
$\D(M)$
in the boson Fock space. Doubtless this representation is irreducible by
 I don't know the precise proof of this (the theorem 8.3 is not quite enough for the proof).

{\bf Problem 8.6.} Representation of $\B(\PSL_2(\R))$
induced from affine isometric action of the group
$\PSL_2(\R)$ is projective. We restrict this representation to the group $\D(M)$ and hence we
 obtain central extensions of
$\D(M)$
enumerated by the points of Teichm\"uller space. But I don't know are these central extensions nontrivial?

{\bf 8.5.Constructions of representations of groups of diffeomorphisms which use $p$-adic numbers.}
 Consider a sphere $S^2$ with volume form $\omega$. Fix points
$m_1,m_2,\dots,m_{l} \in S^2$. Consider the group $\tilde \D(S^2,m)$
of all volume preserving diffeomorphisms of $S^2$ which fix points
$m_1,m_2,\dots,m_{l}$. Let $\D(S^2,{m})$ be the connected component
of unit in the group $\tilde \D(S^2(m)$. The fundamental group
$\pi_1(S^2\setminus\{m_1,m_2,\dots,m_{l}\}$ is a free group
$F_{l-1}$.

Consider the group $ \PSL_2(\K)$ over nonarchimedian local 
field $\K$.
Let us embedd the group $F_{l-1}$ to $ \PSL_2(\K)$ by such a way
that image is a lattice.

The group
$\PSL_2(\K)$
acts on a projective line
$\P\K^1$.
Consider the restriction of this action to
$\Gamma$
and consider the induced action of
$\D(S^2,m)$. 
For this consider the bundle
$\phi : R \rightarrow S^2\setminus\{m_1,\dots,m_l\}$
 associated with principal
$\Gamma$
-bundle
$\H^2 = R \rightarrow S^2\setminus\{m_1,\dots,m_l\}$.
We get a topological space with ergodic action of
$\D(S^2,m)$.
Each
$\D(S^2,m)$
orbit
$\cal O$
is again isomorphic to the Lobachevskii plane and the intersection of
$\cal O$
with each fibre of the bundle is dense in a fibre.

Now we can apply induction procedure.

{\bf 8.6. Once more costruction.}
The factor-group
$\Gamma_m /[\Gamma_m, \Gamma_m]$
is isomorphic to the cyclic group
$\Z^{2m}$.
Consider an action of
$\Z^{2m}$
on
$\C^m$
by translations. It is an affine isometric action and each such action induces the representation of
$\D(M)$
in boson Fock space and the central extension of
$\D(M)$.

\section{On  central extensions of groups of 
automorphisms of bundles}
\setcounter{xxx}{1}
\stepcounter{yyy}

In sections 6-8 we constructed some unitary 
representations and some central extensions of some groups.
It appears that these central extensions exist
in more general situation (when affine isometric 
actions itself  don't exist).

{\bf 9.1. Central extensions of group $\B(\Sp(2n,\R))$.}
Consider space $\R^{2n}$ equipped with skew-symmetric form
defined by the matrix
$\left(\begin{array}{cc}0&1\\-1&0\end{array}\right)$. Denote by
$\Sp(2n,\R)$ the group of operators in $\R^{2n}$which preserve
this form. These operators
$g=\left(\begin{array}{cc}a&b\\c&d\end{array}\right)$
 satisfies the condition
$$\left(\begin{array}{cc}a&b\\c&d\end{array}\right)
\left(\begin{array}{cc}0&1\\-1&0\end{array}\right)
\left(\begin{array}{cc}a&b\\c&d\end{array}\right)^t=
\left(\begin{array}{cc}0&1\\-1&0\end{array}\right)$$

Let $g$ be a element of $\Sp(2n,\R)$. Consider matrix
$$\left(\begin{array}{cc}\Phi&\Psi\\
\overline\Psi&\overline\Phi\end{array}\right):=
\left(\begin{array}{cc}1&i\\i&1\end{array}\right)
\left(\begin{array}{cc}a&b\\c&d\end{array}\right)
\left(\begin{array}{cc}1&i\\i&1\end{array}\right)^{-1}\no$$
We denote $n\times n$-matrix $\Phi=\Phi(g)$ 
as left upper block of the matrix (9.1):
$$\Phi(g)=
\Phi\left\{\left(\begin{array}{cc}a&b\\c&d\end{array}\right)\right\}$$

Consider function $\tau:\Sp(2n,\R)\times\Sp(2n,\R)\rightarrow \R$
given by the formula
$$\tau(g_1,g_2)=
\Im\quad\tr\quad\ln \{\Phi(g_1)^{-1}\Phi(g_1g_2)\Phi(g_2)^{-1}\}$$

{\bf Remark 9.1.} Denote the matrix in curly brackets by  $N$.
Then $\|N-1\|<1$ and hence logarithm of matrix is well-defined. 

Let $(p_1,,h_1), (p_2,h_2)$ be elements of group
 $\B(\Sp(2n,\R))$. We define an expression 
$$c((p_1,,h_1); (p_2,h_2))=
\int_X \tau(h_1(p_2(x)),h_2(x))d\mu(x)$$

{\bf Theorem 9.2.}{\it Expression $\tau$ is 2-cocycle
on the group $\B(\Sp(2n,\R))$.}

{\bf Remark.}
Analogical construction exists for the groups
$\B(U(p,q))$ and $\B(SO^*(2n))$.

{\bf 9.2. Central extensions of groups of symplectomorphisms.}
Let $M$ be a compact $n$-dimensional manifold equipped with symplectic form $\omega$. Let $\SYMP(M,\omega)$ be the group of all diffeomorphisms
of $M$ preserving the form $\omega$. The group 
$\SYMP(M,\omega)$ acts on tangent bundle to $M$ hence we obtain
embedding of $\SYMP(M,\omega)$ to the group $\B(\Sp(2n,\R))$
and this embedding induces some central extension of $\SYMP(M,\omega)$.
We will give description of this extension in coordinates.

Consider some dense open domain $\Xi\subset M$
such that $\Xi$ is symplecticaly equivalent to some domain
$\Omega\subset\R^n$. Now we can consider elements of
$\SYMP$ as transformations of $\Omega$.
Let $g\in\SYMP, x\in\Omega$. Denote by $J(g,x)$ Jacoby
matrix of the map $J$ in the point $x$.
Then 2-cocycle on $\SYMP$ is given by the formula
$$\tau(g_1,g_2)=
\int_\Omega \tau(J(g_1,g_2(x)),J(g_2(x))\omega^n$$

{\bf Proposition 9.3.} {\it Central extension defined by
the 2-cocycle doesn't depend on
coordinate system on M.}

Generally speaking this extension is not trivial,
but I don't know theorem which ca be considered as quite satisfactory.

\section{ Immobile sets and immobile functions}
\setcounter{xxx}{1}
\stepcounter{yyy}

{\bf 10.1. Immobile sets.}
Consider a group $\Gamma$ acting on countable set $M$.
Let $A\subset M$ be a subset. We say that subset $A$
is {\it immobile} if for each $g\in\Gamma$
the set $ A\bigtriangleup gA$ is finite
(where $B\bigtriangleup C:=(B\setminus C)\cup( C\setminus B)$  is
a symmetric sum of the sets).

  Let a group $G$ acts by measure preserving transformations
on space $M$ with infinite measure $\mu$. We say that
subset $A\subset M$ is {\it immobile} if
$$\mu(A\bigtriangleup gA)<\infty$$.

Let $A$ be immobile set and $\mu(A\bigtriangleup  A')<\infty$.
Obviously the set $A'$ also is immobile. We say that immobile sets 
$A$ and $A'$ are {\it equivalent} if
$$ \mu(A\bigtriangleup A')<\infty$$.

Obviously if $A$ and $B$ are immobile sets then $A\cup B$ and
$A\cap B$ are also immobile sets.

{\bf Example 10.1.} Consider a quadratic form $x_1^2+x_2^2-x_3^2$
in $\R^3$. Consider the group ${\rm O}(2,1)$ of all linear transformations preserving this form. Let $G=\SO_0(2,1)$ be the connected component of unit in $SO(2,1)$, see above ss.1.2.
Let $M$ be a one-sheet hyperboloid
$$x_1^2+x_2^2-x_3^2=1$$
and $\mu$ be a $\SO_0(2,1)$-nvariant measure on M. Then
the subset $\{x>0\}$ is $\SO_0(2,1)$-immobile.

{\bf Example 10.2.} Consider the cylindr
$$M:\qquad x_1^2+x_2^2=1$$
in $\R^3$. Consider  function $h(m)=x_3$ on $M$.
Let $G$ be the group of diffeomorphisms $g:M\rightarrow M$ satisfying
condition: $h(m)\rightarrow +\infty$
implies $h(gm)\rightarrow+\infty$. Then the set  
$\{x_3>0\}$ is $G$-immobile.

 Consider a finitely generated discrete group $\Gamma$
acting on countable set $M$. Let $h_1,\dots,h_p$
be generators of the group $\Gamma$. Consider the {\it Cayley graph}
$T(\Gamma,M)$ of our action. Vertices of $ T(\Gamma,M)$ are points
of $M$. Edges have the form $[m,h_jm]$ where
$m\in M$ and $h_j$ are generators.

 {\bf Proposition 10.3} {\it A set $A\subset M$ is immobile iff there exist
only finite number of edges joining $A$ and $M\setminus A$.}

 {\bf Proof}: obvious.

  Let $F_n$ be the free group with $n$ generators $h_1,\dots,h_n$.
Fix a element $v =h_{i_1}^{\varepsilon_1}\dots h_{i_k}^{\varepsilon_k}\in F_n$. Let $X(v)\subset F_n$ be the set of
all irreducible words in $F_n$ having the form
 $$wv=wh_{i_1}^{\varepsilon_1}\dots h_{i_k}^{\varepsilon_k}$$

{\bf Proposition 10.4.} {\it  Consider the left action of the
group $F_n$ on itself. Then the boolean algebra $\Omega$
of all immobile sets is generated by the sets $X(v)$.}

{\bf Proof.} It is a corollary of proposition 10.3.

{\bf 10.2. Immobile functions.} Consider a group $G$ acting on a countable
set $M$. We say that a function $r$ on $M$ is {\it immobile}
if for each $g\in \Gamma$ we have
$$\sum_{m\in M}|r(gm)-r(m)|^2<\infty$$

{\bf Example 10.5.}  Let $A\subset M$ be a immobile set. Then the function
$$\chi(m)=\left\{
\begin{array}{cc} 1,& m\in M
\\2, & m\notin M\end{array}\right.
$$
is immobile.

Let the group $\Gamma$ be finitely generated. Consider
 Cayley graph $T(G,M)$ of our action.

{\bf Proposition 10.6.} {\it Function $r$ on $M$ is immobile iff
$$\sum_{[m,m']\in \Edge(T(G,M))}|r(m)-r(m')|^2<\infty$$
}

{\bf Proof.} Obvious.

Let $r$ be immobile function. Then the formula
$$\Pi(g)f(m)=f(m)+(r(gm)-r(m))$$
defines affine isometric action of the group $\Gamma$ in
the space $l_2(M)$.

{\bf Remark 10.7.} We want to describe the same action by another way.
Fix immobile function $r$. Consider the space of functions
$f$ on $M$ satisfying the condition
$$\sum_m|f(m)-r(m)|^2<\infty$$
Then $H$ has a natural structure of an affine hilbert space.
The group $\Gamma$ acts in $H$ by the transformations
$$\Pi(g)f(m)=f(gm)$$.

{\bf 10.3. On the cohomology group $H^1(\Gamma,l_2(\Gamma)$),} see above  question 0.2.
Consider a discrete group $\Gamma$ and the right regular representation of $\Gamma$
in $l_2(\Gamma)$. Consider all affine actions of $\Gamma$ in
$l_2(\Gamma)$ having the form
$$\Pi(h)f(g)=f(hg)+\gamma_g(h)$$

{\bf Proposition 10.8.}  {\it The cocycle $\gamma_g$ has the form
$$ \gamma_g(h)=r(hg)-r(h)$$
where $r$ is immobile function on $\Gamma$.}

{\bf Proof.} By the condition (0.3)
$$\gamma_{qg}(h)=\gamma_g(qh)+\gamma_q(h)$$
Assume $h=e$. Then
$$\gamma_g(q)=\gamma_{qg}(e)-\gamma(e)$$
and we obtain
$$r(q)=\gamma_q(e)$$.

{\bf The case $\Gamma=F_n$.}\footnote{
When this paper was in preparation I find a preprint
of Bekka and La Valette [BV] which contains more general result 
that a theorem 10.9 formulated below.}
We want to describe the group $H^1(F_n,l_2(F_n))$.
The Cayley graph of action of the group $F_n$ on itself
is the Bruhat--Tits tree ${\cal T}_{2n-1}$.
 Consider the space of functions on the set $\Vert({\cal T}_{2n-1})$
such that
$$\sum_{[p_j,q_j]\in\Edge({\cal T}_{2n-1})}
|f(p_j)-f(q_j)|^2<\infty$$
It is the space $l_2^\pm(\Edge^*)$ defined in ss.1.5. We have to describe quotient space
$$l_2^\pm(\Edge^*)/\Im\bigtriangledown$$
But this space is $\ker \blacktriangledown$. We obtain
following theorem

{\bf Theorem 10.9.}{\it There exist canonical isomorphism
between the space $H^1(F_n, l_2(F_n))$ and the space
$\ker \blacktriangledown=\im {\cal P}$ described in
ss.1.5.}


\end{document}